\documentclass[12pt]{amsart}

\usepackage{amsfonts,amsmath,amssymb}

\title[Classification of Pointed Hopf Algebras]
{On the classification of finite-dimensional pointed Hopf algebras}
\author{Nicol\'{a}s Andruskiewitsch}
\address{Facultad de Matem\'{a}tica, Astronom\'{i}a y f\'{i}sica\\
Universidad Nacional de C\'{o}rdoba \\ CIEM - CONICET,
(5000) Ciudad Universitaria \\
C\'{o}rdoba \\Argentina}
\email{andrus@mate.uncor.edu}
\author{Hans-J\"urgen Schneider}
\address{Mathematisches Institut, Universit\"at M\"unchen, Theresienstr. 39,
D-80333 Munich, Germany}
\email{Hans-Juergen.Schneider@mathematik.uni-muenchen.de}
\thanks{Results of this paper were obtained during a visit of
H.-J. S. at the University of C\'ordoba, partially supported
through a grant of CONICET. The work of N. A. was partially
supported by CONICET, Fundaci\' on Antorchas, Agencia C\'ordoba
Ciencia, ANPCyT, Secyt (UNC) and TWAS (Trieste). This work began
during a stay at the CIRM (Marseille) in the framework of the
program "Research in pairs"; we thank Robert Coquereaux for
providing us with excellent working conditions.}

\newtheorem{Lem}{Lemma}[section]

\newtheorem{Cor}[Lem]{Corollary}
\newtheorem{Thm}[Lem]{Theorem}
\newtheorem{MainThm}[Lem]{Classification Theorem}

\theoremstyle{definition}
\newtheorem{Def}[Lem]{Definition}

\newtheorem{Expl}[Lem]{Example}

\newcommand\pf{\begin{proof}}
\newcommand\epf{\end{proof}}

\newcommand\Aut{\operatorname{Aut}}

\newcommand\id{{\operatorname{id}}}

\newcommand\co{{\operatorname{co}}}

\newcommand\cop{{\operatorname{cop}}}
\newcommand\ad{{\operatorname{ad}}}
\newcommand\het{{\operatorname{ht}}}
\newcommand\gr{{\operatorname{gr}}}
\newcommand\ord{{\operatorname{ord}}}

\newcommand\Isom{{\operatorname{Isom}}}

\renewcommand\o{\otimes}

\newcommand\YDG{^{\Gamma}_{\Gamma}\mathcal{YD}}
\newcommand\YDg{^{\Gamma'}_{\Gamma'}\mathcal{YD}}
\newcommand\YDI{^{\mathbb{Z}[I]}_{\mathbb{Z}[I]}\mathcal{YD}}

\newcommand\ZI{\mathbb{Z}[I]}
\newcommand\D{\mathcal{D}}
\newcommand\X{\mathcal{X}}
\newcommand\Np{\mathbb{N}^p}
\newcommand\ua{\underline{a}}
\newcommand\ub{\underline{b}}
\newcommand\uc{\underline{c}}
\newcommand\ud{\underline{d}}
\newcommand\ue{\underline{e}}
\newcommand\uf{\underline{f}}
\newcommand\ug{\underline{g}}

\newcommand\ur{\underline{r}}
\newcommand\us{\underline{s}}
\newcommand\ut{\underline{t}}

\newcommand\G{\Gamma}
\newcommand\w{\widetilde}
\newcommand\sw[1]{{}_{(#1)}}
\newcommand\swo[1]{{}^{(#1)}}

\numberwithin{equation}{section}

\begin{document}

\maketitle

\begin{abstract}
We classify finite-dimensional complex Hopf algebras $A$ which are
pointed, that is, all of whose irreducible comodules are
one-dimensional, and whose group of group-like elements $G(A)$ is
abelian such that all prime divisors of the order of $G(A)$ are
$>7$. Since these Hopf algebras turn out to be deformations of a
natural class of generalized small quantum groups, our result can
be read as an axiomatic description of generalized small quantum
groups.

\end{abstract}

\section*{Introduction}
One of the very few general classification results for Hopf
algebras is due to Milnor, Moore, Cartier and Kostant around 1963.
It says that any cocommutative Hopf algebra over the complex
numbers is a semidirect product of the universal enveloping
algebra of a Lie algebra and a group algebra.

In the terminology of Sweedler's book \cite{S} from 1969,
cocommutative Hopf algebras (over an algebraically closed field)
are examples of {\em pointed Hopf algebras}, that is, all their
simple subcoalgebras are one-dimensional, or equivalently, all
their simple comodules are one-dimensional. Thus duals of
finite-dimensional pointed Hopf algebras are analogs of basic
algebras in the theory of finite-dimensional algebras. A rich
supply of examples of non-cocommutative pointed Hopf algebra was
only found in the mid-eighties of the last century: The
Drinfeld-Jimbo quantum groups $U_q(\mathfrak{g})$, $\mathfrak{g}$
a semisimple Lie algebra, and their multiparameter versions, as
well as the finite-dimensional small quantum groups introduced by
Lusztig  a bit later are all pointed.

In the present paper we assume that the ground-field $k$ is
algebraically closed of characteristic zero. If $A$ is a Hopf
algebra, we denote the group of group-like elements of $A$ by
$$G(A)= \{g \in A \mid \Delta(g) = g \otimes g, \varepsilon(g) = 1 \}.$$
We classify all finite-dimensional pointed Hopf algebras $A$ over
$k$ with abelian group $G(A)$ such that the prime divisors of the
order of $G(A)$ are $>7$. We describe these Hopf algebras by
generators and relations, and we show that they are the small
quantum groups discovered by Lusztig and variations of them. Thus
our results can be viewed as an axiomatic description of
generalized small quantum groups.

Other types of Hopf algebras occur if prime divisors  $\leq 7$ of
the order of the abelian group $G(A)$ are allowed. However,
extending the methods of this paper it should be possible to
describe their structure by some generalization of Cartan
matrices. Finite-dimensional pointed Hopf algebras with
non-abelian group $G(A)$ seem to be of a very different nature.
Their structure is not understood.

\medskip

We now give a brief overview of the classification program for
finite-dimensional Hopf algebras. We remark that these Hopf
algebras give rise to finite tensor categories in the sense of
\cite{EO} and thus classification results on finite-dimensional
Hopf algebras should have applications in conformal field theory
\cite{Ga}.

The classification splits into several very different parts
according to the behaviour of the coradical. Recall that the
coradical $A_0$ of a Hopf algebra $A$ is the sum of all its simple
subcoalgebras.

\bigskip

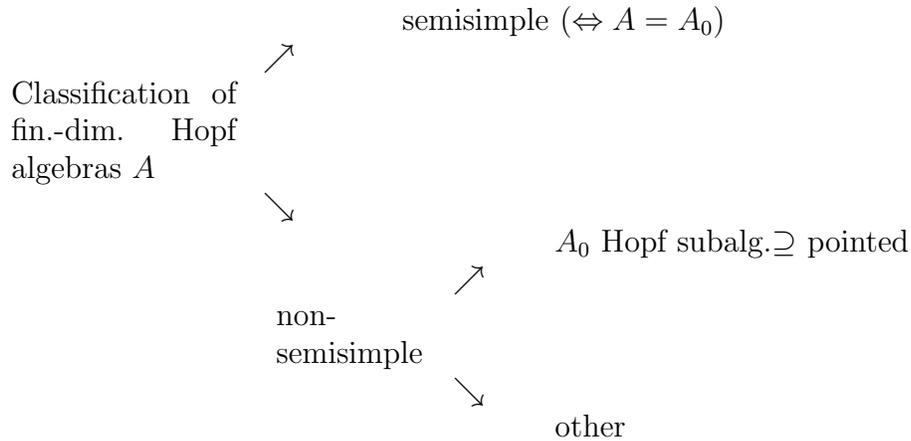
\begin{tabular}{p{3cm}p{1,5cm}p{6cm}}
 && semisimple ($\Leftrightarrow A=A_0$)  \\
& $\nearrow$ & \\
Classification
of fin.-dim. Hopf algebras $A$&&\\
 &$\searrow$ & \\
 &{ \begin{tabular} {p{2cm}p{1cm}p{6cm}}
 && $A_0$  Hopf subalg.$\supseteq$ pointed\\
 &$\nearrow$ & \\
non-semisimple  && \\
& $\searrow$ & \\
 && other
\end{tabular}}
\end{tabular}

\bigskip

If $A=A_0$ then $A$ is semisimple as an algebra. Semisimple Hopf
algebras define examples of fusion categories.  There are various
important results on semisimple Hopf algebras, but there is at
present no general strategy to classify these algebras.

Next assume that $A_0 \neq A$ and that the coradical $A_0$ is a
Hopf subalgebra. If $A$ is pointed then $A_0$ is a Hopf
subalgebra, namely the group algebra $k[G(A)]$. The only general
method for the classification of a class of Hopf algebras is the
Lifting Method, developed in [AS1] and which works for Hopf
algebras whose coradical is a Hopf subalgebra. Since the coradical
is a semisimple Hopf algebra and the classification of semisimple
Hopf algebras is still widely open, it is natural to concentrate
on pointed Hopf algebras. The starting point of this method is an
analog of the Milnor-Moore-Cartier-Kostant decomposition theorem
on the level of the associated graded Hopf algebras. The
enveloping algebra of a Lie algebra is replaced by a braided Hopf
algebra which is generated by primitive elements.

In the case when the coradical $A_0$ of $A$ is not a Hopf
subalgebra very little is known. There are a few results on the
classification of arbitrary Hopf algebras of a given dimension
such as for dimension $p$ or $p^2$ with prime $p$. For their proof
only difficult ad hoc methods are used.

\medskip
To formulate our main result we first describe the data
$\D,\lambda,\mu$ we need to define the Hopf algebras of the class
we are considering. We fix a finite abelian group $\G.$

\subsection*{The datum $\D$} A {\em datum $\D$ of finite Cartan type} for $\G$ ,
$$\mathcal{D} = \mathcal{D}(\Gamma, (g_i)_{1 \leq i \leq \theta}, (\chi_i)_{1 \leq i \leq \theta},
(a_{ij})_{1 \leq i,j \leq \theta}),$$
consists of elements $g_i \in \Gamma, \chi_i \in \widehat{\Gamma}, 1 \leq i \leq \theta,$
and a Cartan matrix $(a_{ij})_{1 \leq i,j \leq \theta}$ of finite type satisfying
\begin{equation}\label{CartantypeIntro}
q_{ij} q_{ji} = q_{ii}^{a_{ij}},\;q_{ii}\neq1, \text{ with } q_{ij} = \chi_j(g_i) \text{ for all } 1 \leq i,j \leq \theta.
\end{equation}
The Cartan condition \eqref{CartantypeIntro} implies  in particular,
\begin{equation}\label{expdiag}
q_{ii}^{a_{ij}} =q_{jj}^{a_{ji}} \text{ for all } 1 \leq i,j \leq \theta.
\end{equation}

The explicit classification of all data of finite Cartan type for
a given finite abelian group $\G$ is a computational problem. But
at least it is a finite problem since the size $\theta$ of the
Cartan matrix is bounded by $2 (\ord(\G))^2$ by \cite[8.1]{AS2},
if $\G$ is an abelian group of odd order. For groups of prime
order, all possibilities for $\D$ are listed in \cite{AS2}.

Let $\Phi$ be the root system of the Cartan matrix $(a_{ij})_{1
\leq i,j \leq \theta},$ $\alpha_1,\dots,\alpha_{\theta}$ a system
of simple roots, and $\mathcal{X}$ the set of connected components
of the Dynkin diagram of $\Phi.$ Let $\Phi_J, J \in \mathcal{X},$
be the root system of the component $J.$ We write $i \sim j,$ if
$\alpha_i$ and $\alpha_j$ are in the same connected component of
the Dynkin diagram of $\Phi.$ For a positive root $\alpha =
\sum_{i=1}^{\theta} n_i \alpha_i, n_i \in
\mathbb{N}=\{0,1,2,\dots\},$ for all $i,$ we define
$$g_{\alpha} = \prod_{i=1}^{\theta} g_i^{n_i}, \chi_{\alpha}
= \prod_{i=1}^{\theta} \chi_i^{n_i}.$$ We assume that the order
of $q_{ii}$ is odd for all $i$, and that the order of $q_{ii}$ is
prime to 3 for all $i$ in a connected component of type $G_2$.
Then it follows from \eqref{expdiag} that the order $N_i$ of
$q_{ii}$ is constant in each connected component $J$, and we
define $N_J = N_i$ for all $i \in J.$

\subsection*{The parameter $\lambda$}
Let $ \lambda = (\lambda_{ij})_{1 \le i < j \le \theta, \, i\not\sim j}$
be a family of elements in $k$ satisfying the following condition
for all $1 \le i < j \le \theta, \, i\not\sim j$: If $g_ig_j =1$ or $\chi_i \chi_j \neq \varepsilon$, then $\lambda_{ij} =0.$

\subsection*{The parameter $\mu$} Let $\mu=(\mu_{\alpha})_{\alpha \in \Phi^+}$ be
a family of elements in $k$ such that for all $\alpha \in
\Phi_J^+,J \in \X,$ if $g_{\alpha}^{N_J} =1$ or
$\chi_{\alpha}^{N_J} \neq \varepsilon,$ then $\mu_{\alpha} =0.$

\medskip
Thus $\lambda$ and $\mu$ are finite families of free parameters in
$k.$ We can normalize $\lambda$ and assume that $\lambda_{ij} =1,$
if $\lambda_{ij} \neq 0.$

\medskip

\subsection*{The Hopf algebra $u(\D,\lambda,\mu)$}
The definition of $u(\D,\lambda,\mu)$ in Section \ref{Sectionu}
can be summarized as follows. In Definition \ref{ualpha} we
associate to any $\mu$ and $\alpha \in \Phi^+$ an element
$u_{\alpha}(\mu)$ in the group algebra $k[\G].$ By construction,
$u_{\alpha}(\mu)$ lies in the augmentation ideal of $k[g_i^{N_i}
\mid 1 \leq i \leq \theta] $. The braided adjoint action
$\ad_c(x_i)$ of $x_i$  is defined in \eqref{braidedfree}, and the
root vectors $x_{\alpha}$ are explained in Section \ref{Cartan}.

\medskip
The Hopf algebra $u(\D,\lambda,\mu)$ is generated as an algebra by
the group $\G,$ that is, by generators of $\G$ satisfying the
relations of the group, and $x_1,\dots,x_\theta,$ with the
relations:
\begin{align*}
&\text{({\em Action of the group}) }& &gx_i g^{-1} = \chi_i(g) x_i,  \text{ for all }
i, \text{ and all } g \in \G,&\\
&\text{({\em Serre relations}) }& &\ad_c(x_i)^{1 - a_{ij}}(x_j) = 0, \text{ for all }  i \neq j, i \sim j,&\\
&\text{({\em Linking relations}) }&& \ad_c(x_i)(x_j) = \lambda_{ij}(1 - g_ig_j), \text{ for all }i<j ,i \nsim j,&\\
&\text{({\em Root vector relations}) }\hspace{-16pt} &&
x_{\alpha}^{N_J} = u_{\alpha}(\mu), \text{ for all }\alpha \in
\Phi_J^+,J \in \mathcal{X}.&
\end{align*}
The coalgebra structure is given by
\begin{align*}
&\Delta (x_i) = g_i \o x_i + x_i \o 1, &&\Delta(g) = g \o g, \text{ for all }1 \leq i \leq \theta,  g \in \G.&
\end{align*}

Now we can formulate our main result.
\begin{MainThm}\label{MainThm}

\noindent (1) Let $\D, \lambda$ and $\mu$ as above. Assume that $q_{ii}$ has
odd order for all $i$ and that the order of $q_{ii}$ is prime to 3
for all $i$ in a connected component of type $G_2.$
Then $u(\D,\lambda,\mu)$ is a pointed Hopf algebra of dimension
$\prod_{J \in  \X} N_J^{|\Phi_J^+|} |\G|$ with group-like elements $G(u(\D,\lambda,\mu)) = \G.$
\medskip

\noindent (2) Let $A$ be a finite-dimensional pointed Hopf algebra
with abelian group $\G = G(A).$ Assume that all prime divisors of
the order of $\Gamma$ are $>7$. Then $A \cong u(\D,\lambda,\mu)$
for some $\D,\lambda,\mu.$
\end{MainThm}

Moreover, in Theorem \ref{iso} we determine all isomorphisms
between the Hopf algebras $u(\D,\lambda,\mu)$.

Part (1) of Theorem \ref{MainThm} is shown in Theorem
\ref{Theoremu(D)}, and part (2) is a special case of Theorem
\ref{Theorempointed}.

In \cite{AS4} we proved the Classification Theorem for groups of
the form $(\mathbb{Z}/(p))^s, s \geq 1,$ where $p$ is a prime
number $>17.$ In this special case, all the elements $\mu$ and
$u_{\alpha}(\mu)$ are zero. In \cite{AS1} we proved part (1) of
Theorem \ref{MainThm} for Dynkin diagrams whose connected
components are of type $A_1$, and in \cite{AS5} for Dynkin
diagrams of type $A_n;$ in \cite{D2} our construction was extended
to Dynkin diagrams whose connected components are of type $A_n$
for various $n$. In \cite{BDR} the Hopf algebra
$u(\D,\lambda,\mu)$ was introduced for type $B_2$.

Our proof of Theorem \ref{MainThm} is based on
\cite{AS1,AS2,AS3,AS4,AS5}, and on previous work on quantum groups
in \cite{dCK,dCP,L1,L2,L3,M1,Ro}, in particular on Lusztig's
theory of the small quantum groups. Another essential ingredient
of our proof are the recent results of Heckenberger on Nichols
algebras of diagonal type in \cite{H1,H2,H3} which use
Kharchenko's theory \cite{K} of PBW-bases in braided Hopf algebras
of diagonal type.

In \cite[1.4]{AS2} we conjectured that any finite-dimensional
pointed Hopf algebra (over an algebraically closed field of
characteristic 0) is generated by group-like and skew-primitive
elements. Our Classification Theorem and Theorem
\ref{Theorempointed} confirm this conjecture for a large class of
Hopf algebras.

Finally we note that the following analog of Cauchy's Theorem from
group theory holds for the Hopf algebras $A=u(\D,\lambda,\mu)$: If
$p$ is a prime divisor of the dimension of $A$, then $A$ contains
a group-like element of order $p$. We conjecture that Cauchy's
Theorem holds for all finite-dimensional pointed Hopf algebras.

\subsection*{Acknowledgement} We thank the referee for very helpful remarks.

\section{Braided Hopf algebras}

\subsection{Yetter-Drinfeld modules over abelian groups and the tensor algebra}\label{tensoralgebra}

Let $\Gamma$ be an abelian group, and $\widehat{\Gamma}$ the
character group of all group homomorphisms from $\Gamma$ to the
multiplicative group $k^{\times}$ of the field $k$.  The braided
category $\YDG$ of (left) Yetter-Drinfeld modules over $\Gamma$
is the category of left $k[\Gamma]$-modules which are
$\Gamma$-graded vector spaces $V = \bigoplus_{g \in \Gamma} V_g$
such that each homogeneous component $V_g$ is stable under the
action of $\Gamma$. Morphisms are $\Gamma$-linear maps $f :
\bigoplus_{g \in \Gamma} V_g \to \bigoplus_{g \in \Gamma} W_g$
with $f(V_g) \subset W_g$ for all $g \in \Gamma$. The
$\Gamma$-grading is equivalent to a left $k[\Gamma]$-comodule
structure $\delta : V \to k[\Gamma] \otimes V$, where $\delta(v) =
g \otimes v$ is equivalent to $v \in V_g$. We use a Sweedler
notation $\delta(v) = v\sw{-1} \otimes v\sw0$ for all $v \in V.$

If $V = \bigoplus_{g \in \Gamma} V_g$ and $W = \bigoplus_{g \in
\Gamma} W_g$ are in $\YDG$, the monoidal structure is given by the
usual tensor product $V \otimes W$ with diagonal $\Gamma$-action $g (v
\otimes w) = gv \otimes gw,\;v \in V,w \in W$, and
$\Gamma$-grading $(V \otimes W)_g = \bigoplus_{ab=g} V_a\otimes
W_b$ for all $g \in \Gamma.$ The braiding in $\YDG$ is the
isomorphism
$$c = c_{V,W} : V \otimes W \to W \otimes V$$
defined by $c(v \otimes w) = g \cdot w \otimes v$ for all $g \in
\Gamma, v \in V_g, \text{ and } w \in W.$ Thus each
Yetter-Drinfeld module $V$ defines a braided vector space
$(V,c_{V,V}).$

If $\chi$ is a character of $\Gamma$  and $V$ a left
$\Gamma$-module, we define
$$V^{\chi} := \{ v \in V \mid  g \cdot v = \chi(g) v \text{ for all } g \in \Gamma\}.$$
Let $\theta \geq 1$ be a natural number, $g_1,\dots,g_{\theta} \in
\Gamma,$ and $\chi_1,\dots,\chi_{\theta} \in \widehat{\Gamma}.$
Let $V$ be a vector space with basis $x_1,\dots, x_{\theta}.$ $V$
is an object in $\YDG$ by defining $x_i \in V^{\chi_i}_{g_i}$ for
all $i$. Thus each $x_i$ has degree $g_i$, and the group $\Gamma$
acts on $x_i$ via the character $\chi_i$. We define
$$q_{ij} := \chi_j(g_i) \text{  for all } 1 \leq i,j \leq \theta.$$
The braiding on $V$ is determined by the matrix $(q_{ij})$ since
$$c(x_i \o x_j) = q_{ij} x_j \o x_i \text{ for all } 1 \leq i,j \leq \theta.$$
We will identify the tensor algebra $T(V)$ with the free
associative algebra $k\langle x_1,\dots,x_\theta\rangle$. It is an
algebra in $\YDG,$ where a monomial
$$x = x_{i_1}x_{i_2} \cdots x_{i_n}, 1 \leq i_1, \dots , i_n \leq \theta,$$
has $\Gamma$-degree $g_{i_1}g_{i_2} \cdots g_{i_n}$ and the action
of $g \in \Gamma$ on $x$ is given by $g \rightharpoonup x
=\chi_{i_1}\chi_{i_2} \cdots \chi_{i_n}(g)x.$ $T(V)$ is a braided
Hopf algebra in $\YDG$ with comultiplication
$$\Delta_{T(V)} : T(V) \to T(V) \underline{ \otimes} T(V), \;
x_i \mapsto x_i \otimes 1 + 1 \otimes x_i,\; 1 \leq i \leq \theta.$$
Here we write $T(V) \underline{ \otimes} T(V)$ to indicate the
braided algebra structure on the vector space $T(V) \otimes T(V)$,
that is
$$(x \otimes y)(x' \otimes y') = x (g \rightharpoonup x') \otimes yy',$$ for all
$x,x',y,y' \in T(V)$ and $y \in T(V)_g, g \in \Gamma.$

\medskip
Let $I=\{1,2,\dots,\theta\}$, and $\mathbb{Z}[I]$ the free abelian
group of rank $\theta$ with basis $\alpha_1,\dots,
\alpha_{\theta}.$ Given the matrix $(q_{ij})$, we define the
bilinear map
\begin{equation}\label{qalpha}
\ZI \times \ZI \to k^{\times},\; (\alpha,\beta) \mapsto q_{\alpha \beta},
\text{ by } q_{\alpha_i \alpha_j} = q_{ij}, 1 \leq i,j \leq \theta.
\end{equation}
We consider $V$ as a Yetter-Drinfeld module over $\mathbb{Z}[I]$
with $x_i \in V_{\alpha_i}^{\psi_i}$ for all $1 \leq i \leq
\theta$, where $\psi_j$ is the character of $\mathbb{Z}[I]$ with
$$\psi_j(\alpha_i) = q_{ij} \text{  for all } 1 \leq i,j \leq \theta.$$
Thus $T(V) = k\langle x_1,\dots,x_\theta \rangle$ is also a
braided Hopf algebra in $\YDI$. The $\ZI$-degree of a monomial $x
= x_{i_1}x_{i_2} \cdots x_{i_n}, 1 \leq i_1, \dots , i_n \leq
\theta,$ is $\sum_{i=1}^{\theta} n_i \alpha_i$, where for all $i$,
$n_i$ is the number of occurences of $i$ in the sequence
$(i_1,i_2,\dots,i_n)$. It follows for the action of $\mathbb{Z}[I]$ on homogeneous components that
\begin{equation}\label{actionhom}
\alpha \rightharpoonup x = q_{\alpha \beta}x \text{ for all } \alpha,\beta \in \mathbb{Z}[I], x \in T(V)_{\beta}.
\end{equation}
The braiding on $T(V)$ as a Yetter-Drinfeld
module over $\Gamma$ or $\ZI$ is in both cases given by
\begin{equation}\label{freebraiding}
c(x \o y) = q_{\alpha \beta}y \o x, \text{ where } x \in T(V)_{\alpha}, y \in T(V)_{\beta}, \alpha, \beta \in \ZI.
\end{equation}
The comultiplication of $T(V)$ as a braided Hopf algebra in $\YDG$
only depends on the matrix $(q_{ij})$, hence it coincides with the
comultiplication of $T(V)$ as a coalgebra in $\YDI$. In
particular, the comultiplication of $T(V)$ is
$\mathbb{Z}[I]$-graded.

\subsection{Bosonization and twisting}\label{twisting}

Let $R$ be a braided Hopf algebra in $\YDG$. We will use a
Sweedler notation for the comultiplication
$$\Delta_R : R \to R \otimes R, \; \Delta_R(r) = r \swo1 \otimes r \swo2.$$
For Hopf algebras $A$ in the usual sense, we always use the Sweedler notation
$$\Delta : A \to A \otimes A, \; \Delta(a) = a \sw1 \otimes a \sw2.$$
Then the smash product $A = R \# k[\Gamma]$ is a Hopf algebra in
the usual sense (the bosonization of $R$). As vector spaces, $R \#
k[\Gamma] = R \otimes k[\Gamma]$. Multiplication and
comultiplication are defined by
\begin{equation}\label{smashproduct}
(r \# g)(s \# h) = r (g\cdot s) \# gh, \;
\Delta(r \# g) = r \swo1 \# r \swo2 \sw{-1}g \otimes r \swo2 \sw0 \# g.
\end{equation}
Then the maps
$$\iota : k[\Gamma] \to R \# k[\Gamma],  \text{ and } \pi : R \# k[\Gamma] \to k[\Gamma]$$
with $\iota(g) = 1 \# g \text{ and } \pi(r \# g) = \varepsilon(r)g$ for all $r
\in R$, $g \in \Gamma$ are Hopf algebra maps with $\pi \iota =
\id.$

For simplicity we will often write $rg$ instead of $r \# g$ in $R \# k[\G]$ for $r \in R,g \in \G$. Thus $\Delta_{R \# k[\G]}(r)=r \swo1 r\swo2 \sw{-1} \o r \swo2 \sw0$.

Conversely, if $A$ is a Hopf algebra in the usual sense with Hopf
algebra maps $\iota : k[\Gamma] \to A \text{ and } \pi : A \to
k[\Gamma]$ such that $\pi \iota = \id,$ then
\begin{equation}\label{defR}
R= \{a \in A \mid (\id \o \pi)\Delta(a) = a \otimes 1\}
\end{equation} is a
braided Hopf algebra in $\YDG$ in the following way. As an
algebra, $R$ is a subalgebra of $A$. The $k[\Gamma]$-coaction,
$\Gamma$-action and comultiplication of $R$ are defined by
\begin{equation}\label{Radfordcoaction}
\delta(r) = \pi(r_{(1)}) \otimes r_{(2)}, \quad g \rightharpoonup
r = \iota(g) r \iota(g^{-1})
\end{equation}
and
\begin{equation}\label{Radfordcomult}
 \Delta_R(r) = \vartheta(r \sw1) \otimes r \sw2.
\end{equation}
Here, $\Delta_A(r) = r \sw1 \otimes r \sw2,$ and $\vartheta$ is the map
\begin{equation}
\vartheta : A \to R,\; \vartheta(r) = r\sw1 \iota (S(\pi(r\sw2))),
\end{equation}
where $S$ is the antipode of $A$. Then
\begin{equation}\label{Radfordiso}
R \# k[\G] \to A, \; r \# g \mapsto r \iota(g),\; r \in R, g \in \G,
\end{equation}
is an isomorphism of Hopf algebras.

\medskip
We recall the notion of {\em twisting} the algebra structure of an
arbitrary Hopf algebra $A$, see for example \cite[10.2.3]{KS}. Let
$\sigma : A \otimes A \to k$ be a convolution invertible linear
map, and a normalized 2-cocycle, that is, for all $x,y,z \in A,$
\begin{equation}\label{cocycle}
\sigma(x\sw1,y\sw1) \sigma(x\sw2 y\sw2, z) = \sigma(y\sw1,z\sw1) \sigma(x,y\sw2 z\sw2),
\end{equation}
and $\sigma(x,1) = \varepsilon(x) = \sigma(1,x).$ The Hopf algebra
$A_{\sigma}$ with twisted algebra structure is equal to $A$ as a
coalgebra, and has multiplication $\cdot_{\sigma}$ with
\begin{equation}\label{twistedmultiplication}
x \cdot_{\sigma}y = \sigma(x\sw1,y\sw1) x\sw2 y\sw2 \sigma^{-1}(x\sw3,y\sw3) \text{ for all } x,y \in A.
\end{equation}
In the situation $A = R \# k[\Gamma]$ above, let $\sigma : \Gamma
\times \Gamma \to k^{\times}$ be a normalized 2-cocycle of the
group $\Gamma$. Then $\sigma$ extends to a 2-cocycle of the group
algebra $k[\Gamma]$ and it defines a normalized and invertible
2-cocycle
$$\sigma_{\pi} = \sigma(\pi \otimes \pi)$$
of the Hopf
algebra $A$. Since $k[\Gamma]$ is cocommutative, $\iota$ and $\pi$
are Hopf algebra maps
$$\iota : k[\Gamma] \to A_{\sigma_{\pi}} \text{ and } \pi : A_{\sigma_{\pi}} \to k[\Gamma].$$
Hence the coinvariant elements
$$R_{\sigma} = \{ a \in A_{\sigma_{\pi}} \mid (\id \o \pi)\Delta(a) = a \otimes 1\}$$
form a braided Hopf algebra in $\YDG$. As a $k[\G]$-comodule,
$R_{\sigma}$ coincides with $R$,
but $R_{\sigma}$ and $R$ have different action, multiplication and
comultiplication. To simplify the formulas, we will treat $\iota$
as an inclusion map. Also we denote the new action by
$\rightharpoonup_{\sigma}$.

In any braided Hopf algebra $R$ with multiplication $m$ and
braiding $c : R \otimes R \to R \otimes R$ we define the {\em
braided commutator} of elements $x,y \in R$ by
\begin{equation}\label{commutator}
[x,y]_c = xy - mc(x \otimes y).
\end{equation}
If $x \in R$ is a primitive element, then
\begin{equation}\label{braidedadjoint}
(\ad_cx)(y) = [x,y]_c
\end{equation}
denotes the {\em braided adjoint action} of $x$ on $R$. For
example, in the situation of the free algebra in Section
\ref{tensoralgebra} with braiding \eqref{freebraiding}, we have
for all  $x_i$ and $y = x_{j_1} \cdots x_{j_n},$
\begin{equation}\label{braidedfree}
(\ad_cx_i)(y) = x_iy - q_{ij_1} \cdots q_{ij_n}yx_i.
\end{equation}
To formulate the next lemma we need one more notation. If
$V$ is a left $C$-comodule over a coalgebra $C$, then $V$ is a
right module over the dual algebra $C^*$ by $v \leftharpoonup p =
p(v\sw{-1}) v\sw{0}$ for all $v \in V, p\in C^*.$ In particular,
if $R$ is a braided Hopf algebra in $\YDG$, then the
$k[\Gamma]$-coaction defines a left $k[\Gamma]\otimes
k[\Gamma]$-comodule structure on $R \otimes R$, hence a right
$(k[\Gamma]\otimes k[\Gamma])^*$-module structure on $R \otimes R$
denoted by $\leftharpoonup.$

\begin{Lem}\label{braidedtwist}
Let $\Gamma$ be an abelian group, $\sigma : \Gamma \times \Gamma
\to k^{\times}$  a normalized 2-cocycle, $R$ a braided Hopf
algebra in $\YDG$, $g,h \in \Gamma,x \in R_g,y\in R_h$, and $r
\in R.$ Then
\begin{enumerate}
\item $x \cdot_{\sigma} y = \sigma(g,h) xy.$\label{twistedmult}
\item $\Delta_{R_{\sigma}}(r) = \Delta_R(r) \leftharpoonup \sigma^{-1}.$\label{twistedcomult}
\item\label{twistedcomm} If $y \in R_h^{\eta}$ for some character $\eta \in \widehat{\Gamma}$,
then
\begin{equation}\label{twistedaction}
g \rightharpoonup_{\sigma} y = \sigma(g,h) \sigma^{-1}(h,g)
\eta(g) y,
\end{equation}
and hence $[x,y]_{c_{\sigma}} = \sigma(g,h)[x,y]_{c}.$
\end{enumerate}
\end{Lem}
\pf Note that for all homogeneous elements $z \in R_s,s \in \G,$
$$(\pi\otimes
\id\otimes \pi) \Delta^2(z) = s\otimes z \otimes 1,$$
because of
\eqref{defR} and \eqref{Radfordcoaction}. This implies \eqref{twistedmult}, and in \eqref{twistedcomm} we obtain
$$g\cdot_{\sigma}y
= \sigma(g,h) gy \text{ and }y\cdot_{\sigma}g = \sigma(h,g) yg.$$
Thus
$$g\cdot_{\sigma}y = \sigma(g,h) gy = \sigma(g,h)\eta(g) yg =
\sigma(g,h)\eta(g)\sigma(h,g)^{-1}y\cdot_{\sigma}g$$
and \eqref{twistedaction} follows. In turn, \eqref{twistedaction}
implies the last assertion in \eqref{twistedcomm}.

To prove
\eqref{twistedcomult}, using the cocommutativity of the group
algebra we compute
\begin{align*}
\Delta_{R_{\sigma}}(r) &= r\sw1 \cdot_{\sigma} S(\pi(r\sw2)) \otimes r\sw3\\
&=\sigma(\pi(r\sw1),S(\pi(r\sw5))) \vartheta(r\sw2)\sigma^{-1}(\pi(r\sw3),S(\pi(r\sw4))) \otimes r\sw6.
\end{align*}
On the other hand, $\Delta_R(r) = r\sw1 S\pi(r\sw2) \otimes
r\sw3,$ hence

\noindent $r\swo1\sw{-1} \otimes r\swo2\sw{-1} \otimes r\swo1\sw0
\otimes r\swo2\sw0 =\pi(r\sw1 S(r\sw3)) \otimes \pi(r\sw4) \otimes
\vartheta(r\sw2) \otimes r\sw5,$ and  $\Delta_R(r) \leftharpoonup
\sigma^{-1} = \sigma^{-1}(\pi(r\sw1 S(r\sw3)), \pi(r\sw4))
\vartheta(r\sw2) \otimes r\sw5.$ Hence the claim follows from the
equality
$$\sigma(a,S(b\sw3)) \sigma^{-1}(b\sw1,S(b\sw2)) = \sigma^{-1}(a S(b\sw1),b\sw2)$$
for all $a,b \in k[\Gamma].$ It is enough to check this equation
for elements $a,b \in \Gamma$. Then the equality follows from the
group cocycle condition, which implies $\sigma(g, g^{-1}) =
\sigma(g^{-1}, g)$ for all $g\in \Gamma$.
\epf
Part (3) of the previous lemma extends \cite[Lemma 3.15]{AS4} and \cite[Lemma 2.12]{AS5}.

We now apply the twisting procedure to the braided Hopf algebra
$T(V) \in {\YDI}.$

\begin{Lem}\label{mapphi}
Let $\theta\geq 1$, and let $(q_{ij})_{1\leq i,j\leq
\theta},(q'_{ij})_{1\leq i,j\leq \theta}$ be matrices with
coefficients in $k$. As in Section \ref{tensoralgebra} let $V \text{ and } V' \in {\YDI}$ with basis
$x_1,\dots,x_{\theta}$ and $x'_1,\dots,x'_{\theta}$ respectively, where $x_i \in V_{\alpha_i}^{\psi_i},x'_i \in {V'}_{\alpha_i}^{{\psi'}_i}$ with
$\psi_j(\alpha_i) = q_{ij}, \psi'_j(\alpha_i) = q'_{ij}$ for all $i,j$. Then $T(V)$ and
$T(V')$ are braided Hopf algebras in $\YDI$ as in Section
\ref{tensoralgebra}. Assume
\begin{equation}\label{equivbraiding}
q_{ij}q_{ji}=q'_{ij}q'_{ji}, \text{ and } q_{ii} = q'_{ii} \text{ for all } 1 \leq i,j \leq \theta.
\end{equation}
Then there is a 2-cocycle $\sigma : \mathbb{Z}[I] \times \mathbb{Z}[I] \to k^{\times}$ with
\begin{equation}\label{equivcocycle}
\sigma(\alpha,\beta) \sigma^{-1}(\beta,\alpha) = q_{\alpha\beta}
q'^{-1}_{\alpha\beta} \text{ for all } \alpha, \beta \in \ZI,
\end{equation}
and a $k$-linear isomorphism $\varphi : T(V) \to T(V')$  with
$\varphi(x_i) = x'_i$ for all $i$ and such that for all $\alpha,
\beta \in \mathbb{Z}[I], x \in T(V)_{\alpha}, y \in T(V)_{\beta}$
and $z \in T(V)$
\begin{enumerate}
\item $\varphi(xy) = \sigma(\alpha,\beta) \varphi(x) \varphi(y).$\label{equivmult}
\item $\Delta_{T(V')}(\varphi(z)) = (\varphi\otimes \varphi)(\Delta_{T(V)}(z))
\leftharpoonup \sigma.$\label{equivcomult}
\item $\varphi([x,y]_c) = \sigma(\alpha,\beta)[\varphi(x),\varphi(y)]_{c'}.$\label{equivcommutator}
\end{enumerate}
\end{Lem}
\pf Define $\sigma$ as the bilinear map with
$\sigma(\alpha_i,\alpha_j) = q_{ij}q'^{-1}_{ij}$ if $i\leq j$, and
$\sigma(\alpha_i,\alpha_j) = 1$ if $i>j$ (see \cite[Prop.
3.9]{AS5}).

Let $\varphi : T(V) \to T(V')_{\sigma}$ be the algebra map with
$\varphi(x_i) = x'_i$ for all $i$. Then $\varphi$ is bijective
since it follows from Lemma \ref{braidedtwist} \eqref{twistedmult}
and the bilinearity of $\sigma$ that for all monomials $x =
x_{i_1} x_{i_2} \cdots x_{i_n}$ of length $n \geq 1$ with $x' =
x'_{i_1} x'_{i_2} \cdots x'_{i_n},$
$$\varphi(x) = \prod_{r<s} \sigma(\alpha_{i_r}, \alpha_{i_s}) x'.$$
In particular, $\varphi$ is $\mathbb{Z}[I]$-graded.

To see that
$\varphi$ is $\mathbb{Z}[I]$-linear, let $\alpha,\beta \in
\mathbb{Z}[I]$ and $x \in T(V)_{\beta}.$ By \eqref{actionhom} and Lemma
\ref{braidedtwist} \eqref{twistedcomm},
$$\alpha \rightharpoonup x = q_{\alpha\beta} x,\text{ and } \alpha
\rightharpoonup_{\sigma} \varphi(x) = \sigma(\alpha,\beta)
\sigma^{-1}(\beta,\alpha) q'_{\alpha\beta} \varphi(x),$$ and
$\varphi(\alpha \rightharpoonup x) = \alpha \rightharpoonup_{\sigma} \varphi(x)$
follows by \eqref{equivcocycle}.

Since the elements $x_i$ and
$x'_i$ are primitive we now see that $\varphi$ is an isomorphism of braided Hopf algebras. Then
the claim follows from Lemma \ref{braidedtwist}. \epf

\section{Serre relations and root vectors}\label{Serre}

\subsection{Datum of finite Cartan type and root vectors}\label{Cartan}

\begin{Def}
A {\em datum of Cartan type}
$$\mathcal{D} = \mathcal{D}(\Gamma, (g_i)_{1 \leq i \leq \theta},
(\chi_i)_{1 \leq i \leq \theta}, (a_{ij})_{1 \leq i,j \leq
\theta})$$ consists of an abelian group $\Gamma$, elements $g_i
\in \Gamma, \chi_i \in \widehat{\Gamma}, 1 \leq i \leq \theta,$
and a generalized Cartan matrix $(a_{ij})$ of size $\theta$ satisfying
\begin{equation}\label{Cartantype}
q_{ij} q_{ji} = q_{ii}^{a_{ij}},\;q_{ii}\neq 1, \text{ with }
q_{ij} = \chi_j(g_i) \text{ for all } 1 \leq i,j \leq \theta.
\end{equation}
We call $\theta$ the {\em rank} of $\D$. A datum $\mathcal{D}$ of
Cartan type will be called of {\em finite Cartan type} if $(a_{ij})$ is
of finite type.
\end{Def}
\begin{Expl}\label{ELusztig}
A Cartan datum $(I, \cdot)$  in the sense of Lusztig
\cite[1.1.1]{L3} defines a datum of Cartan type for the free
abelian group $ZI$ with $g_i = \alpha_i, \chi_i = \psi_i, 1 \leq i
\leq \theta,$ as in Section \ref{tensoralgebra}, where
$$ q_{ij} = v^{d_ia_{ij}}, d_i = \frac{i \cdot i}{2}, a_{ij} = 2 \frac{i \cdot j}{i \cdot i}
\text{ for all } 1 \leq i,j \leq \theta.$$
\end{Expl}
In Example \ref{ELusztig}, $d_ia_{ij} = i \cdot j$ is the
symmetrized Cartan matrix, and $q_{ij} = q_{ji}$ for all $1 \leq
i,j \leq \theta.$ In general, the matrix $(q_{ij})$ of a datum of
Cartan type is not symmetric, but by Lemma \ref{mapphi} we can
reduce to the symmetric case by twisting.

\medskip

We fix a finite abelian group $\Gamma$ and a datum
$$\mathcal{D} = \mathcal{D}(\Gamma, (g_i)_{1 \leq i \leq \theta},
(\chi_i)_{1 \leq i \leq \theta}, (a_{ij})_{1 \leq i,j \leq
\theta})$$ of finite Cartan type. The Weyl group $W \subset
\Aut(\ZI)$ of $(a_{ij})$ is generated by the reflections $s_i :
\ZI \to \ZI$ with $s_i(\alpha_j) = \alpha_j - a_{ij} \alpha_i$ for
all $i,j.$ The root system is $\Phi = \cup_{i = 1}^{\theta}
W(\alpha_i),$ and
$$\Phi^+ = \Big\{\alpha \in \Phi \mid \alpha = \sum_{i=1}^{\theta} n_i
\alpha_i, n_i \geq 0 \text{ for all } 1 \leq i \leq \theta\Big\}$$
denotes the set of positive roots with respect to the basis of
simple roots $\alpha_1,\dots, \alpha_{\theta}.$ Let $p$ be the
number of positive roots.

For $\alpha = \sum_{i=1}^{\theta} n_i \alpha_i \in \ZI, n_i \in \mathbb{Z} \text{ for all } i$, we define
\begin{equation}\label{galpha}
g_{\alpha} = g_1^{n_1} g_2^{n_2} \cdots g_{\theta}^{n_{\theta}} \text{ and }
\chi_{\alpha} = \chi_1^{n_1} \chi_2^{n_2} \cdots \chi_{\theta}^{n_{\theta}}.
\end{equation}
Hence for all $ \alpha,\beta \in \ZI$,
\begin{equation}\label{alphabeta}
q_{\alpha \beta} = \chi_{\beta}(g_{\alpha}),
\end{equation}
where $q_{\alpha \beta}$ is given by \eqref{qalpha}.

In this section, we assume that the Dynkin diagram of $(a_{ij})$
is {\em connected}. In this case we say that $\D$ is connected.
We assume for all $1 \leq i \leq \theta,$
\begin{align}
&q_{ii} \text{ has odd order} , \text{ and }&\label{orderodd}\\
&\text{the order of } q_{ii} \text{ is prime to 3, if $(a_{ij})$
is of type } G_2.&\label{orderG2}
\end{align}
Then it follows from Lemma \ref{q} that the elements $q_{ii}$ have
the same order in $k^{\times}$. We define
\begin{equation}\label{orderN}
N = \text{ order of } q_{ii}, 1 \leq i \leq \theta.
\end{equation}

\begin{Lem}\label{q}
Let $\D$ be a connected datum of finite Cartan type with Cartan
matrix $(a_{ij})$ and assume \eqref{orderodd} and \eqref{orderG2}. Then there are integers $d_i \in \{1,2,3\}, 1 \leq i
\leq \theta$, and $q \in k$ such that for all $1 \leq i,j \leq
\theta$,
$$q_{ii} = q^{2d_i},\quad d_ia_{ij}=d_ja_{ji},$$
and the order of $q$ is odd, and if the Cartan matrix of $\D$ is of type $G_2$, then the order of $q$ is prime to 3.
\end{Lem}
\pf By \eqref{Cartantype}, $q_{ii}^{a_{ij}} = q_{jj}^{a_{ji}}$ for
all $1 \leq i,j \leq \theta$. It then follows from the list of
Cartan matrices with connected Dynkin diagram that in each case
there is an index $h$ such that $q_{ii} =(q_{hh})^{d_i}$ for all
$1 \leq i \leq \theta$, where the $d_i \in \{1,2,3\}$ symmetrize
$(a_{ij})$. Indeed this is obvious in the simply laced case, and
in the notation of \cite{B} with $l=\theta$ take $h=1$ for $B_l$,
$h=l$ for $C_l$, and $h=2$ for $F_4$ and $G_2$. By \eqref{orderodd} we let $q$ be
a square root of $q_{hh}$ of odd order. Then by \eqref{orderG2} the order of $q$ is prime to 3 if the Cartan matrix of $\D$ is of type $G_2$.
\epf

We fix a reduced decomposition of the longest element
$$w_{0}= s_{i_1} s_{i_2} \cdots s_{i_p}$$
of $W$ in terms of the simple reflections.
Then
$$\beta_{l} = s_{i_1} \cdots s_{i_{l-1}}(\alpha_{i_l}), 1 \leq l \leq p,$$
is a convex ordering of the positive roots.

\begin{Def}\label{DefR(D)}
Let $V = V(\mathcal{D})$ be a vector space with basis
$x_1,\dots,x_\theta$, and let $V \in {\YDG}$ by $x_i \in
V_{g_i}^{\chi_i}$ for all $1 \leq i \leq \theta.$ Then $T(V)$ is a
braided Hopf algebra in $\YDG$ as in Section \ref{tensoralgebra}.
Let
$$R(\D) = T(V)/((\ad_cx_i)^{1-a_{ij}}(x_j) \mid 1\leq i\neq j \leq \theta)$$
be the quotient Hopf algebra in $\YDG.$
\end{Def}

It is well-known that the elements $(\ad_cx_i)^{1-a_{ij}}(x_j),
1\leq i\neq j \leq \theta$, are primitive in the free algebra
$T(V)$ (see for example \cite[A.1]{AS2}), hence they generate a
Hopf ideal. By abuse of language, we denote the images of the
elements $x_i$ in $R(\D)$ again by $x_i$.

\medskip
In the situation of Example \ref{ELusztig}, Lusztig \cite{L2}
defined root vectors $x_\alpha$ in $R(\D)= U^+$ for each positive
root $\alpha$ using the convex ordering of the positive roots. As
noted in \cite{AS4}, these root vectors can be seen to be iterated
braided commutators of the elements $x_1,\dots,x_{\theta}$ with
respect to the braiding given by the matrix $(v^{d_ia_{ij}})$.
This follows for example from the inductive definition of the root
vectors in \cite{Ri}.

In the case of our general braiding given by $(q_{ij})$ we define
root vectors $x_{\alpha} \in R(\D)$ for each $\alpha \in \Phi^+$
by the same iterated braided commutator of the elements
$x_1,\dots,x_{\theta}$ as in Lusztig's case but with respect to
the general braiding.

\begin{Def}\label{DefK(D)}
Let $K(\D)$ be the subalgebra of $R(\D)$ generated by the elements $x_{\alpha}^N, \alpha \in \Phi^+$.
\end{Def}

\begin{Thm}\label{R(D)}
Let $\D$ be a connected datum of finite Cartan type, and assume \eqref{orderodd}, \eqref{orderG2}.
\begin{enumerate}
\item The elements
$$x_{\beta_1}^{a_1}x_{\beta_2}^{a_2} \cdots x_{\beta_p}^{a_p}, a_1,a_2,\dots,a_p \geq 0,$$
form a basis of $R(\D)$.
\item $K(\D)$ is a braided Hopf subalgebra of $R(\D).$
\item For all $\alpha, \beta \in \Phi^+, [x_{\alpha}, x_{\beta}^N]_c = 0,$ that is,
$$x_{\alpha} x_{\beta}^N =
q_{\alpha \beta}^N x_{\beta}^N x_{\alpha}.$$
\end{enumerate}
\end{Thm}
\pf (a) In the situation of Example \ref{ELusztig}, the elements
in (1) form Lusztig's PBW-basis of $U^+$ over
$\mathbb{Z}[v,v^{-1}]$ by \cite[5.7]{L2}.

(b) Now we assume that the braiding has the form $(q_{ij} =
q^{d_ia_{ij}})$, where $(d_ia_{ij})$ is the symmetrized Cartan
matrix, and $q$ is a non-zero element in $k$ of odd order, and whose order is prime to 3 if the Dynkin diagram of $(a_{ij})$ is $G_2$. Then
(1) follows from Lusztig's result by extension of scalars, and (2)
is shown in \cite[19.1]{dCP} (for another proof see
\cite[3.1]{M2}). The algebra $K(\D)$ is commutative since it is a
subalgebra of the commutative algebra $Z_0$ of \cite[19.1]{dCP}.
This proves (3) since  $q^N = 1$, hence
$\chi_{\beta}^N(g_{\alpha}) = 1$.

(c) In the situation of a general braiding matrix $(q_{ij})_{1\leq
i,j \leq \theta}$ assumed in the theorem, we apply Lemma \ref{q}
and define a matrix $(q'_{ij})_{1\leq i,j \leq \theta}$ by
$q'_{ij} = q^{d_ia_{ij}}$ for all $i,j$. Then $q_{ij}q_{ji} =
q'_{ij}q'_{ji}$, and $q_{ii}=q'_{ii}$ for all $1 \leq i,j \leq
\theta$. Thus by part (b) of the proof, (1),(2) and (3) hold for
the braiding $(q'_{ij})$, and hence by Lemma \ref{mapphi} for
$(q_{ij})$. \epf

\subsection{The Hopf algebra $K(\D) \# k[\Gamma]$}\label{SectionK(D)}

We assume the situation of Section \ref{Cartan}. By Theorem
\ref{R(D)} (2), $K(\D)$ is a braided Hopf algebra in $\YDG$, and
the smash product $K(\D) \# k[\Gamma]$ is a Hopf algebra in the
usual sense. We want to describe all Hopf algebra maps
$$K(\D) \# k[\Gamma] \to k[\Gamma]$$
which are the identity on the group algebra $k[\Gamma].$

\begin{Def}\label{Defza}
For any $1 \leq l \leq p$ and $a = (a_1,a_2, \dots, a_p) \in \mathbb{N}^p$ we define
\begin{align*}
&&h_l &= g_{\beta_l}^N, &\\
&&\eta_l &= \chi_{\beta_l}^N, &\\
&&z_l &= x_{\beta_l}^N,&\\
&&z^a &= z_{1}^{a_1}z_{2}^{a_2} \cdots z_{_p}^{a_p} \in K(\D),&\\
&&h^a &= h_1^{a_1}h_{2}^{a_2} \cdots h_{p}^{a_p} \in \Gamma,&\\
&&\eta^a &= \eta_{1}^{a_1}\eta_{_2}^{a_2} \cdots \eta_{p}^{a_p} \in \widehat{\Gamma},&\\
&&\underline{a} &= a_1 \beta_1 + a_2 \beta_2 + \cdots + a_p \beta_p \in \ZI.&
\end{align*}
For $\alpha = \sum_{i=1}^{\theta} n_i \alpha_i \in \ZI, n_i \in
\mathbb{Z} \text{ for all } i,$ we call $\het(\alpha) =
\sum_{i=1}^{\theta} n_i$ the {\em height} of $\alpha.$ Let $e_l =
(\delta_{kl})_{1 \leq k \leq p} \in \Np$, where $\delta_{kl} = 1$
if $k=l$ and $\delta_{kl} =0$ if $k \neq l.$
\end{Def}

Note that for all $a,b,c \in \mathbb{N}^p$,
\begin{equation}
h^a = h^b h^c,\;\eta^a = \eta^b  \eta^c, \text{ if } \underline{a} = \underline{b} + \underline{c},
\end{equation}
\begin{equation}
\het(\underline{b}) < \het(\underline{a}), \text{ if }
\underline{a} = \underline{b} + \underline{c} \text{ and }c \neq
0.
\end{equation}

\bigskip
As explained in Section \ref{tensoralgebra}, we view $T(V)$ as a
braided Hopf algebra in $\YDI$. Then the quotient Hopf algebra
$R(\D)$ and its Hopf subalgebra $K(\D)$ are braided Hopf algebras
in $\YDI.$ In particular, the comultiplication $\Delta_{K(\D)} :
K(\D) \to K(\D) \o K(\D)$ is $\ZI$-graded. By construction, for
any $\alpha \in \Phi^+,$ the root vector $x_{\alpha}$ in $R(\D)$
is $\ZI$-homogeneous of $\ZI$-degree $\alpha.$ Thus $x_{\alpha}
\in R(\D)_{g_{\alpha}}^{\chi_{\alpha}},$ and for all $a \in
\mathbb{N}^p,$ $z^a$ has $\ZI$-degree $N \underline{a},$ and
\begin{equation}\label{degreez}
z^a \in K(\D)_{h^a}^{\eta^a}.
\end{equation}
By Theorem \ref{R(D)} the elements $z^ag$
with $a \in \mathbb{N}^p, g \in \Gamma,$ form a basis of $K(\D) \#
k[\Gamma],$ and it follows that for all $a, b = (b_i),c = (c_i)\in
\mathbb{N}^p,$
\begin{equation}\label{multiplicationK(D)}
z^b z^c = \gamma_{b,c} z^{b + c}, \text{ where } \gamma_{b,c} = \prod_{k>l} \eta_l(h_k)^{b_kc_l},
\end{equation}
\begin{equation}\label{groupactionK(D)}
h^az^b = \eta^b(h^a) z^b h^a \text{ in } R \# k[\Gamma].
\end{equation}

\begin{Lem}\label{formulat}
For any $0 \neq a \in \mathbb{N}^p$ there are uniquely determined
scalars $t_{b,c}^a \in k,0 \neq  b,c \in \mathbb{N}^p,$ such that
\begin{equation}\label{comultiplicationformula}
\Delta_{K(\D)}(z^a) = z^a \o 1 + 1 \o z^a + \sum_{b,c \neq 0,
\underline{b} + \underline{c} = \underline{a}} t_{b,c}^a \,z^b \o
z^c.
\end{equation}
\end{Lem}
\pf Since $\Delta_{K(\D)}$ is $\ZI$-graded, $\Delta_{K(\D)}(z^a)$
is a linear combination of elements $z^b \o z^c$ where
$\underline{b} + \underline{c} = \underline{a}$. Hence
$$\Delta_{K(\D)}(z^a) = x \o 1 + 1 \o y + \sum_{b,c \neq 0,
\underline{b} + \underline{c} = \underline{a}} t_{b,c}^a\, z^b \o
z^c,$$ where $x,y$ are elements in $K(\D)$. By applying the
augmentation $\varepsilon$ it follows that $x = y = z^a.$ \epf

We now define recursively a family of elements $u^a$ in $k[\G]$
depending on parameters $\mu_a$ which behave like the elements
$z^a$ with respect to comultiplication.

\begin{Lem}\label{inductionstep1}
Let $n \geq 1$. Let $(\mu_b)_{ 0 \neq b \in \Np, \het(\ub) <n}$ be a family of elements in $k$, and  let $(u^b)_{ 0 \neq b \in \Np, \het(\ub) <n}$ be a family of elements in $k[\G]$. Assume for all $0 \neq b \in \Np, \het(\ub) <n,$ that
\begin{equation}\label{ub}
u^b = \mu_b(1 - h^b) + \sum_{d,e \neq 0, \ud + \ue = \ub} t_{d,e}^b\, \mu_d u^e,
\end{equation}
\begin{equation}\label{Deltaub}
\Delta(u^b) = h^b \o u^b + u^b \o 1 + \sum_{d,e \neq 0, \ud + \ue = \ub} t_{d,e}^b\, u^d h^e \o u^e.
\end{equation}
Let $a \in \Np$ with $\het(\ua) = n,$ and $u^a \in k[\Gamma]$.
Then the following statements are equivalent:
\begin{equation}\label{ua}
u^a = \mu_a(1 - h^a) + \sum_{b,c \neq 0, \ub + \uc = \ua} t_{b,c}^a \,\mu_b u^c \text{ for some } \mu_a \in k.
\end{equation}
\begin{equation}\label{Deltaua}
\Delta(u^a) = h^a \o u^a + u^a \o 1 + \sum_{b,c \neq 0, \ub + \uc = \ua} t_{b,c}^a \,u^b h^c \o u^c.
\end{equation}
\end{Lem}
\pf
If $n=1$, the equivalence between \eqref{ua} and \eqref{Deltaua}
is well-known and easy to see. The point of the Lemma is the inductive
construction of the $u^a$'s.
Let $$v_a = u^a - \sum_{b,c \neq 0, \ub + \uc = \ua}
t_{b,c}^a\, \mu_b u^c.$$ Then $u^a$ can be written as in
\eqref{ua} if and only if $\Delta(v_a) = h^a \o v_a + v_a \o 1.$
Hence it is enough to prove that
$$\Delta(v_a) - h^a \o v_a - v_a \o 1 = \Delta(u^a) - h^a \o u^a -
u^a \o 1 - \sum_{b,c \neq 0, \ub + \uc = \ua} t_{b,c}^a \,u^b h^c
\o u^c.$$ We compute
\begin{align}
\Delta(v_a)&- h^a \o v_a - v_a \o 1 = \notag\\
&=\Delta(u^a) - \sum_{b,c \neq 0, \ub + \uc = \ua} t_{b,c}^a\, \mu_b \Delta(u^c) - h^a \o v_a - v_a \o 1  \notag\\
&=\Delta(u^a) - h^a \o u^a - u^a \o 1 + \sum_{b,c \neq 0, \ub + \uc = \ua} t_{b,c}^a\, \mu_b(h^a \o u^c - h^c \o u^c)\notag \\
& - \sum_{\substack{b,c,f,g \neq 0\\ \ub + \uc =\ua, \uf + \ug = \uc}} t_{b,c}^a\, t_{f,g}^c \,\mu_b u^f h^g \o u^g,\notag
\end{align}
using the definition of $v_a$ in the first equation, and the
formula for $\Delta(u^c)$ from \eqref{Deltaub} in the second
equation. Note that the term
$$\sum_{b,c \neq 0, \ub + \uc = \ua} t_{b,c}^a \,\mu_b u^c \o 1$$
cancels. Hence we have to show that
\begin{align}
&\sum_{\substack{b,c,f,g \neq 0\\ \ub + \uc =\ua, \uf + \ug = \uc}} t_{b,c}^a\, t_{f,g}^c\, \mu_b u^f h^g \o u^g =\notag\\
&= \sum_{b,c \neq 0, \ub + \uc = \ua} t_{b,c}^a (\mu_bh^a \o u^c - \mu_bh^c \o u^c + u^b  h^c \o u^c).\notag
\end{align}
Since for all $b,c \neq 0, \ub + \uc = \ua$, we have $h^a = h^b
h^c,$ it follows that
$$\mu_bh^a \o u^c - \mu_bh^c \o u^c + u^b h^c \o u^c = (\mu_b(h^b - 1) + u^b)h^c \o u^c.$$
Using the formula for $u^b$ from \eqref{ub}, we finally have to prove
$$\sum_{\substack{b,c,f,g \neq 0\\ \ub + \uc =\ua, \uf + \ug = \uc}} t_{b,c}^a \,t_{f,g}^c\, \mu_b u^f h^g \o u^g = \sum_{\substack{b,c,d,e \neq 0\\ \ub + \uc = \ua, \ud + \ue = \ub}}t_{b,c}^a\, t_{d,e}^b\, \mu_d u^e h^c \o u^c.$$
This last equality follows from the coassociativity of $K(\D)$. Indeed, from
$$(\id \o \Delta_{K(\D)})\Delta_{K(\D)}(z^a) = (\Delta_{K(\D)}\o \id)\Delta_{K(\D)}(z^a)$$
we obtain with \eqref{comultiplicationformula} after cancelling
several terms
$$\sum_{\substack{b,c,f,g \neq 0\\ \ub + \uc =\ua, \uf + \ug = \uc}} t_{b,c}^a\, t_{f,g}^c \,z^b \o z^f \o z^g  = \sum_{\substack{b,c,d,e \neq 0\\ \ub + \uc = \ua, \ud + \ue = \ub}}t_{b,c}^a\, t_{d,e}^b\, z^d \o z^e \o z^c.$$
Thus mapping $z^r \o z^s \o z^t, r,s,t \neq 0, \het(\ur),
\het(\us),\het(\ut) <n,$ onto $\mu_r u^s h^t \o u^t$ proves the
claim. Here we are using that the elements $z^a$ are linearly
independent by Theorem \ref{R(D)}. \epf

Let $K(\D) \# k[\Gamma]$ be the Hopf algebra corresponding to the
braided Hopf algebra $K(\D)$ by \eqref{smashproduct}.
 Thus by definition and Lemma \ref{formulat}, for all $0 \neq a \in \Np,$
 \begin{equation}\label{smashcomult}
\Delta_{K(\D) \# k[\Gamma]}(z^a) = h^a \o z^a + z^a \o 1 + \sum_{b,c \neq 0, \ub + \uc = \ua} t_{b,c}^a\, z^b h^c \o z^c.
\end{equation}

For all $n \geq 0$, let $K(\D)_n$ be the vector subspace spanned
by all elements $z^a, a \in \Np, \het(\ua) \leq n.$ Then $K(\D)_n \#
k[\Gamma] \subset K(\D) \# k[\Gamma]$ is a subcoalgebra.

\bigbreak In the next Lemma we describe all coalgebra maps
$$\varphi : K(\D)_n \# k[\Gamma] \to k[\Gamma] \text{ with } \varphi | \Gamma = \id.$$
Note that such a coalgebra map is given by a family of elements
$\varphi(z^a) = :u^a, 0\neq a \in \Np, \het(\ua) \leq n,$ such
that  \eqref{Deltaua} holds for all $0\neq a, \het(\ua) \leq n.$
It follows by induction on $\het(\ua)$ from Lemma
\ref{inductionstep1} with \eqref{ua} that $\varepsilon(u^a) = 0$
for all $a.$

\begin{Lem}\label{coalgebramaps}
Let $n \geq 1.$

{\rm (1)} Let $(\mu_a)_{0 \neq a \in \Np, \het(\ua) \leq  n}$ be a
family of elements in $k$ such that for all $a$, if $h^a =1$, then
$\mu_a =0$. Define the family $(u^a)_{0 \neq a \in \Np, \het(\ua)
\leq n}$ by induction on $\het(\ua)$ by \eqref{ua}. Then the map
$\varphi : K(\D)_n \# k[\Gamma] \to k[\Gamma]$ given by $\varphi | \Gamma = \id$, 
$$  \varphi(z^ag) =
u^ag, 0 \neq a \in \Np, \het(\ua) \leq  n, g \in \Gamma,$$ is a coalgebra
map.

\medskip
{\rm (2)} The map defined in {\rm(1)} from the set of all $(\mu_a)_{0 \neq a
\in \Np, \het(\ua) \leq n}$ such that for all $a$, if $h^a =1$,
then $\mu_a =0$, to the set of all coalgebra maps $\varphi$ with
$\varphi | \Gamma = \id$ is bijective.
\end{Lem}
\pf This follows from Lemma \ref{inductionstep1} by induction on
$\het(\ua).$ Note that the coefficient $\mu_a$ in \eqref{ua} is
uniquely determined if we define $\mu_a = 0$ if $h^a = 1.$ \epf

\begin{Def}\label{partialHopf}
Let $n \geq 1.$ A coalgebra map
$$\varphi : K(\D)_n \# k[\Gamma]
\to k[\Gamma]
\text{ with }\varphi |\Gamma = \id$$ is called a {\em
partial Hopf algebra map}, if for all $x,y \in K(\D)_n \#
k[\Gamma]$ with $xy \in K(\D)_n \# k[\Gamma]$, we have
$\varphi(xy) = \varphi(x) \varphi(y).$
\end{Def}

\begin{Lem}\label{partial}
Let $n \geq 1,$ and $\varphi : K(\D)_n \# k[\Gamma] \to k[\Gamma]$
a coalgebra map, $(\mu_a)_{0 \neq a \in \Np, \het(\ua) \leq n}$
the family of scalars corresponding to $\varphi$ by Lemma
\ref{coalgebramaps}, and $u^a = \varphi(z^a)$ for all $a \in \Np$
with $\het(\ua) \leq n.$ Then the following are equivalent:
\begin{enumerate}
\item $\varphi$ is a partial Hopf algebra map.\label{partialmap}
\item For all $0 \neq a= (a_1,\dots,a_p) \in \Np$ with $\het(\ua) \leq n,$
\begin{enumerate}
\item $u^a = \prod_{a_l>0} u_l^{a_l},$ where for all $1 \leq l \leq p,u_l = u^{e_l},$ if $a_l >0,$
\item if $\eta^a \neq \varepsilon,$ then $\mu_a = 0,$ and $u^a = 0.$
\end{enumerate}
\item \begin{enumerate}
\item As {\upshape (2) (a)}.
\item For all $1 \leq l \leq p$ with $\het(\underline{e_l}) \leq n,$ if
$\eta_l \neq \varepsilon$, then $u^{e_l} = 0.$
\end{enumerate}
\end{enumerate}
\end{Lem}
\pf $(1) \Rightarrow (2)$: If $\varphi$ is a partial Hopf algebra
map, then (a) follows immediately, and to prove (b), let $0 \neq a
\in \Np, \het(\ua) \leq n,$ and $g \in \Gamma,$ with $\eta^a \neq
\varepsilon.$  Then
$$\varphi(gz^a) = \eta^a(g) u^ag = u^ag,$$
since $gz^a = \eta^a(g) z^ag$ by \eqref{groupactionK(D)}.  Thus
$u^a =0$, and it follows by induction on $\het(\ua)$ from
\eqref{ua} that $\mu_a =0,$ since for all $ 0 \neq b,c \in \Np$
with $\het(\ub) + \het(\uc) = \het(\ua)$, $\eta^b \neq
\varepsilon,$ or $\eta^c \neq \varepsilon.$

\medskip
$(2) \Rightarrow (3)$ is trivial. $(3) \Rightarrow (1)$: The
coalgebra map $\varphi$ is a partial Hopf algebra map if and only
if for all $b,c \in \Np$ with $\het(\ub) + \het(\uc) \leq n,$ and
$g,h \in \Gamma,$
$$\varphi(z^bg z^ch) = u^bg u^c h.$$
By \eqref{multiplicationK(D)} and \eqref{groupactionK(D)},
$z^bg z^ch = \eta^c(g) \gamma_{b,c} z^{b+c} gh.$ Thus (1) is equivalent to
\begin{equation}\label{equivalentpartial}
\eta^c(g) \gamma_{b,c} u^{b+c} = u^b u^c \text{ for all } b,c \in \Np, \het(\ub) + \het(\uc) \leq n, g\in \Gamma.
\end{equation}
Let $b,c \in \Np, \het(\ub) + \het(\uc) \leq n, g\in \Gamma.$ By (a),
$$u^{b+c} = u^b u^c = \prod_{b_l + c_l >0} u_l^{b_l+c_l}.$$
To prove \eqref{equivalentpartial} assume that $u^b u^c \neq 0.$
Then $u_l \neq 0$ for all $l$ with $c_l >0.$ Hence by (b), $\eta_l
= \varepsilon$ for all $l$ with $c_l >0,$ and
$\eta^c(g)=1,\gamma_{b,c}=1$. \epf

To formulate the main result of this section, we define $M(\D)$ as
the set of all families $(\mu_l)_{1 \leq l \leq p}$ of elements in
$k$ satisfying the following condition for all $ 1 \leq l \leq p:$
If $h_l=1$ or $\eta_l \neq \varepsilon,$ then $\mu_l =0.$

\begin{Thm}\label{mainconstruction}
{\rm(1)} Let $\mu = (\mu_l)_{1 \leq l \leq p} \in M(\D)$. Then there is
exactly one Hopf algebra  map
$$ \varphi_{\mu} : K(\D) \# k[\Gamma] \to k[\Gamma], \;\varphi|\Gamma = \id$$
such that the family $(\mu_a)_{0 \neq a \in \Np}$ associated to
$\varphi_{\mu}$ by Lemma \ref{coalgebramaps} satisfies  $\mu_{e_l}
= \mu_l$ for all $1 \leq l \leq p.$

\medskip
{\rm(2)} The map $\mu \mapsto \varphi_{\mu}$ defined in {\rm(1)} from
$M(\D)$ to the set of all Hopf algebra homomorphisms $\varphi :
K(\D) \# k[\Gamma] \to k[\Gamma]$ with $ \varphi|\Gamma = \id$ is
bijective.
\end{Thm}
\pf (1) We proceed by induction on $n$ to construct partial Hopf
algebra maps on $K(\D)_n \# k[\G],$ the case $n=0$ being trivial.
We assume that we are given a partial Hopf algebra map
$$\varphi : K(\D)_{n-1} \# k[\Gamma] \to k[\Gamma], \; n \geq 1,$$
such that $\mu_{e_l} = \mu_l$ for all $1\leq l \leq p$ with
$\het(\underline{e_l}) \leq n-1.$ Here $(\mu_a)_{0\neq a \in \Np,
\het(\ua) \leq n-1}$ is the family of scalars associated to
$\varphi$ by Lemma \ref{coalgebramaps}. We define $u^b =
\varphi(z^b)$ for all $0 \neq b, \het(\ub) \leq n-1$. It is enough
to show that there is exactly one partial Hopf algebra map
$$\psi: K(\D)_{n} \# k[\Gamma] \to k[\Gamma]$$
extending $\varphi,$ and such that $\mu_{e_l} = \mu_l$ for all $l$ with $\het(\underline{e_l}) \leq n.$

Let $a \in \Np$ with $\het(\ua) =n.$ To define $\psi(z^a)=: u^a$
we distinguish two cases.

If $a = e_l$ for some $1 \leq l \leq p$, we define
\begin{equation}\label{defua}
u^a = \mu_l(1 - h^a) + \sum_{b,c \neq 0, \ub + \uc = \ua} t_{b,c}^a \mu_b u^c.
\end{equation}
Then \eqref{Deltaua} holds by Lemma \ref{inductionstep1}.

If $a = (a_1, \dots, a_l,0,\dots,0), a_l  \geq 1, 1\leq l \leq p,$
and $a \neq e_l$, then $a = r + s$, where $0\neq r, s = e_l.$ We
define $u^a = u^r u^s.$ To see that $u^a$ satisfies
\eqref{Deltaua}, using \eqref{smashcomult} we write
$$\Delta(z^c) = h^c \o z^c + z^c \o 1 + T(c), \text{ for all } 0 \neq c \in \Np.$$
Since $z^r z^s = z^a$ because of \eqref{multiplicationK(D)} (note
that $\gamma_{r,s} = 1$ in this case) we see that
$\Delta(z^r)\Delta(z^s) = h^a \o z^a + z^a \o 1 + T(r,s),$ where
\begin{multline*}
T(r,s) = h^r z^s \o z^r + z^rh^s \o z^s\\
+ (h^r \o z^r + z^r \o 1)T(s) +T(r)(h^s \o z^s + z^s \o 1) +
T(r)T(s),
\end{multline*}
and $T(r,s) = T(a).$ Since $\varphi$ on $K(\D)_{n-1} \# k[\Gamma]$
is a coalgebra map,
$$\Delta(u^c) = h^c \o u^c + u^c \o 1 +(\varphi \o \varphi)(T(c)),$$
for all $0 \neq c \in \Np$ with $\het(\uc) \leq n-1.$
In particular,
$$\Delta(u^r) \Delta(u^s) = h^a \o u^a + u^a \o 1 + (\varphi \o \varphi)(T(r,s)).$$
Thus $\Delta(u^a) = h^a \o u^a + u^a \o 1 + (\varphi \o
\varphi)(T(a)),$ that is, $u^a$ satisfies \eqref{Deltaua}.

Thus the extension of $\varphi$ defined by $\psi(z^ag) = u^ag$ for
all $g \in \Gamma, a \in \Np, \het(\ua) = n$ is a coalgebra map.

\medskip To prove that the extension $\psi$ is a partial Hopf algebra map,
we check condition (3) in Lemma \ref{partial}. Since the
restriction of $\psi$ to $K(\D)_{n-1} \# k[\Gamma]$ is a partial
Hopf algebra map, (3) (a) is satisfied. To prove (3)(b), let $1
\leq l \leq p$ with $\het(\underline{e_l}) =n,$  $a = e_l,$ and
assume $\eta_l \neq \varepsilon.$ Then for all $0 \neq b,c \in
\Np$ with $\ub + \uc = \ua$, we have $\eta^b \neq \varepsilon$ or
$\eta^c \neq \varepsilon.$ Since $\varphi$ is a Hopf algebra map,
it follows from Lemma \ref{partial} that $\mu_b =0$ or $u^c =0.$
By assumption, $\mu_l = 0.$ Hence by \eqref{defua}, $u^a = 0$.

This proves (1) since the uniqueness of the extension follows from
Lemma \ref{inductionstep1} and Lemma \ref{coalgebramaps}.

\medskip (2) By Lemma \ref{coalgebramaps}, the map $\mu \mapsto
\varphi_{\mu}$ is injective. To prove surjectivity, let $\varphi :
K(\D) \# k[\Gamma] \to k[\Gamma]$ be a Hopf algebra map with
$\varphi | \Gamma = \id.$ By Lemma \ref{coalgebramaps}, $\varphi$
is defined by a family $(\mu_a)_{0 \neq a \in \Np}$ of scalars. By
(1), $\varphi$ is determined by the values $\mu_{e_l}, 1 \leq l
\leq p.$ \epf

\begin{Def}\label{ualpha}
For any $\mu \in M(\D)$ and $1 \leq l \leq p$, let $\varphi_{\mu}$
be the Hopf algebra map defined in Theorem \ref{mainconstruction},
and
$$u_l(\mu) = \varphi_{\mu}(z_l) \in  k[\Gamma].$$
If $\alpha$ is a positive root in $\Phi^+$ with $\alpha =
\beta_l$, we define $u_{\alpha}(\mu) = u_l(\mu).$
\end{Def}

Note that by \eqref{ua}, each $u_{\alpha}(\mu)$ lies in the
augmentation ideal of $k[g_{i}^N \mid  1 \leq i \leq \theta].$

\section{Linking}\label{Sectionlinking}

\subsection{Notations}\label{Sectionnotations}

In this Section we fix a finite abelian group $\Gamma$, and a
datum $\D = \D(\Gamma, (g_i)_{1 \leq i \leq \theta}, (\chi_i)_{1
\leq i \leq \theta}, (a_{ij})_{1\leq i,j \leq \theta})$ of finite
Cartan type. We follow the notations of the previous Section, in
particular,
$$q_{ij} = \chi_j(g_i)\text{  for all }i,j.$$

For all $ 1 \leq i,j \leq \theta$ we write $i \sim j$ if $i$ and
$j$ are in the same connected component of the Dynkin diagram of
$(a_{ij})$. Let $\mathcal{X} = \{I_1,\dots,I_t\}$ be the set of
connected components of $I = \{1,2, \dots, \theta\}$. We assume for all $1 \leq i \leq \theta$
\begin{align}\label{orderqij}
&q_{ii} \text{ has odd order}, \text{ and}&\\
&\text{the order of } q_{ii} \text{ is prime to 3, if $i$ lies in a component }G_2.\label{orderNJ}
\end{align}
For all $J \in \X,$ let $N_J$ be the common order of $q_{ii}, i \in J.$

As in Section \ref{SectionK(D)}, for all $J \in \X$ we choose a
reduced decomposition of the longest element $w_{0,J}$ of the Weyl
group $W_J$ of the root system $\Phi_J$ of $(a_{ij})_{i,j \in J}$.
Then for all $J,K \in \X,$ $w_{0,J}$ and $w_{0,K}$ commute in the
Weyl group $W$ of the root system $\Phi$ of $(a_{ij})_{1\leq i,j
\leq \theta}$, and
$$w_0 = w_{0,I_1} w_{0,I_2} \cdots w_{0,I_t}$$
gives a reduced representation of the longest element of $W$. For
all $J \in \X,$ let $p_J$ be the number of positive roots in
$\Phi_J^+$, and
\begin{equation*}\label{orderPhil}
\Phi_{J}^+ = \{ \beta_{J,1}, \dots, \beta_{J,p_{J}}\}
\end{equation*}
the corresponding convex ordering. Then
\begin{equation*}\label{orderPhi}
\Phi^+ = \{ \beta_{I_1,1},\dots, \beta_{I_1,p_{I_1}},\dots,\beta_{I_t,1},\dots, \beta_{I_t,p_{I_t}}\}
\end{equation*}
is the convex ordering corresponding to the reduced representation
of $w_0 = w_{0,I_1} w_{0,I_2} \cdots w_{0,I_t}$. We also write
\begin{equation*}\label{orderPhialso}
\Phi^+ = \{\beta_1,\dots,\beta_p\},\; p = \sum_{J \in \X} p_{J},
\end{equation*}
for this ordering.

In Section \ref{Cartan} we have defined root vectors $x_{\alpha}$
in the free algebra $k\langle x_1,\dots,x_{\theta}\rangle$ for
each positive root in $\Phi_J^+ \subset \Phi, J \in \X.$

\medskip
We recall a notion from \cite{AS4}.

\begin{Def}\label{Deflinking}
A family $ \lambda = (\lambda_{ij})_{1 \le i <j \le \theta, \,
i\not\sim j}$ of elements in $k$ is called a {\em family of
linking parameters for $\D$} if the following condition is
satisfied for all $1 \le i , j \le \theta, \, i\not\sim j$: 
\begin{equation}
\text{If }g_ig_j =1 \text{ or }\chi_i \chi_j \neq \varepsilon, \text{ then }
\lambda_{ij} =0.\label{l1}
\end{equation}
Vertices $1\leq i,j \leq \theta$ are called
{\em linkable} if $i \not\sim j,$ $g_ig_j \neq 1$ and $\chi_i
\chi_j = \varepsilon$.
\end{Def}
It is useful to formally extend the notion of linking parameters by 
\begin{equation}\label{l2}
\lambda_{ji} = -q_{ji} \lambda_{ij}\text{ for all }1 \le i <j \le \theta, \,
i\not\sim j.
\end{equation}

Assume in addition that $\ord(q_{ii}) >3$ for all $i$. We remark that any vertex $i$ is linkable to at most one vertex $j$, and if $i,j$
are linkable, then $q_{ii}= q_{jj}^{-1}$ \cite[Section 5.1]{AS4}.

\medskip

The free algebra $k\langle x_1,\dots,x_{\theta}\rangle$ is a
braided Hopf algebra in $\YDG$ as explained in Section
\ref{tensoralgebra}. Then $k\langle x_1,\dots,x_{\theta}\rangle \#
k[\Gamma]$ is a Hopf algebra as in \ref{twisting}. 

\subsection{The Hopf algebra $U(\D, \lambda)$}\label{SectionU(D)}

We assume the situation of Section \ref{Sectionnotations}.
\begin{Def}\label{DefU}
Let $ \lambda = (\lambda_{ij})_{1 \le i , j \le \theta, \,
i\not\sim j}$ be a family of linking parameters for $\D$. Let
$U(\D,\lambda)$ be the quotient Hopf algebra of  \newline
$k\langle x_1,\dots,x_{\theta}\rangle \# k[\Gamma]$ modulo the
ideal generated by
\begin{align}
&\ad_c(x_i)^{1-a_{ij}}(x_j), \text{ for all } 1\leq i,j \leq \theta, i \sim j, i \neq j,&\label{Serrerelations}\\
&x_i x_j - q_{ij} x_j x_i - \lambda_{ij}(1 - g_i g_j), \text{ for all } 1 \le i < j \le \theta, \, i\not\sim j.&\label{linkingrelations}
\end{align}
\end{Def}
In \eqref{linkingrelations} we can add the redundant elements
$$x_j x_i - q_{ji} x_i x_j - \lambda_{ji}(1 - g_j g_i), 1 \le i < j \le \theta, \, i\not\sim j.$$ Thus the definition of $U(\D,\lambda)$ does not depend on the ordering of the index set.

We denote the images of $x_i$ and $g\in \Gamma$ in $U(\D,\lambda)$
again by $x_i$ and $g$. The elements in \eqref{Serrerelations} and
\eqref{linkingrelations} are skew-primitive. Hence
$U(\D,\lambda)$ is a Hopf algebra with
$$\Delta(x_i) = g_i \o x_i + x_i \o 1,\; 1 \leq i \leq \theta.$$

In part (1) of the next theorem we adapt the method of proof of \cite[Section 5.3]{AS4} to find a basis of the infinite-dimensional Hopf algebra $U(\D,\lambda)$ in terms of the root vectors. In part (2) we prove a crucial skew-commutativity relation for the root vectors.

\begin{Thm}\label{TheoremU(D)}
Let $\Gamma$ be a finite abelian group, and $\D$ a datum of finite
Cartan type satisfying \eqref{orderqij} and \eqref{orderNJ}. Let
$\lambda$ be a family of linking parameters for $\D$. Then
\begin{enumerate}
\item The elements
$$x_{\beta_1}^{a_1} x_{\beta_2}^{a_2} \cdots x_{\beta_p}^{a_p}g,\; a_1, a_2,\dots, a_p \geq 0, g \in \Gamma,$$
form a basis of the vector space $U(\D,\lambda).$
\item Let $J \in \X$ and $\alpha \in \Phi^+,  \beta \in \Phi_{J}^+.$
Then $[x_{\alpha},x_{\beta}^{N_J}]_c = 0,$ that is,
$$x_{\alpha} x_{\beta}^{N_J} = q_{\alpha \beta}^{N_J} x_{\beta}^{N_J} x_{\alpha}.$$
\end{enumerate}
\end{Thm}
\pf We proceed by induction on the number $t$ of connected components.

If $I$ is connected, (1) and (2) follow from Theorem \ref{R(D)}.

If $t >1$, let $I_1 =
\{1,2,\dots,\widetilde{\theta}\},\; 1 \leq \widetilde{\theta} <
\theta.$ For all $1 \leq i \leq \widetilde{\theta},$ let $l_i$ be
the least common multiple of the orders of $g_i$ and $\chi_i$, $1
\leq i \leq \widetilde{\theta}.$ Let $\widetilde{\Gamma} = \langle
\widetilde{g}_1,\dots,\widetilde{g}_{\widetilde{\theta}}\mid \widetilde{g}_i \widetilde{g}_j = \widetilde{g}_j \widetilde{g}_i,
\widetilde{g}_i^{l_i}=1 \text{ for all } 1 \leq i,j \leq \widetilde{\theta}\rangle$, and for all $1
\leq j \leq \widetilde{\theta}$ let $\widetilde{\chi}_j$ be the character of
$\widetilde{\Gamma}$ with $\widetilde{\chi}_j(\widetilde{g}_i) = \chi_j(g_i)$ for all $1\leq i
\leq \widetilde{\theta}.$ Then we define
\begin{align*}
\D_1 &= \D(\widetilde{\Gamma},(\widetilde{g}_i)_{1 \leq i \leq
\widetilde{\theta}},(\widetilde{\chi}_i)_{1 \leq i \leq
\widetilde{\theta}},(a_{ij})_{1 \leq i,j \leq
\widetilde{\theta}}),\\
\D_2 &= \D(\Gamma,
(g_i)_{\widetilde{\theta} <  i \leq \theta},
(\chi_i)_{\widetilde{\theta} <  i \leq \theta},
(a_{ij})_{\widetilde{\theta} < i,j \leq \theta}),
\end{align*}
and $\lambda_2
=(\lambda_{ij})_{\widetilde{\theta} < i < j \leq \theta, i \nsim
j}.$ Let $U = U(\D_1)$ (with empty family of linking
parameters) with generators $x_1,\dots,x_{\widetilde{\theta}}$,
and $\widetilde{g} \in \widetilde{\Gamma}$, and $A = U(\D_2,\lambda_2)$ with
generators $y_{\widetilde{\theta} +1}, \dots, y_{\theta},$ and $g
\in \Gamma.$

It is shown in \cite[Lemma 5.19]{AS4} that there are algebra maps
$\gamma_i$, $(\varepsilon,\gamma_i)-$derivations $\delta_i$ and a
Hopf algebra map $\varphi,$
$$\gamma_i : A \to k,\; \delta_i : A \to k,\; \varphi : U \to
(A^0)^{\cop},\; 1 \leq i \leq \widetilde{\theta},$$ such that for
all $1\leq i \leq \widetilde{\theta}< j \leq \theta,$
\begin{alignat*}{2}
\gamma_i | \Gamma &= \chi_i,\quad& \gamma_i(y_j) &= 0,\\
\delta_i |\Gamma &=0, &\delta_i(y_j) &= \lambda_{ji},\\
\varphi(\widetilde{g}_i) &= \gamma_i,& \varphi(x_i) &= \delta_i.
\end{alignat*}
Then $\sigma : U \o A \o U \o A \to U \o A,$ defined for all $u,v \in U, a,b \in A$ by
\begin{enumerate}
 \item [] $\sigma(u \o a,v \o b) = \varepsilon(u) \tau(v,a) \varepsilon(b),\;\tau(v,a)=\varphi(v)(a), $
\end{enumerate}
is a 2-cocycle on the tensor product Hopf algebra of $U$ and $A,$
and $(U \o A)_{\sigma}$ is the Hopf algebra with twisted
multiplication defined in \eqref{twistedmultiplication}.
Multiplication in $(U \o A)_{\sigma}$ is given for all $u,v \in U,
a,b \in A$ by
\begin{equation}\label{multiplicationrule}
(u \o a) \cdot_{\sigma} (v \o b) = u\tau(v\sw1, a\sw1) v\sw2 \o a\sw2 \tau^{-1}(v\sw3, a\sw3) b,
\end{equation}
with $\tau^{-1}(u,a) = \varphi(u)(S^{-1}(a)).$

The group-like elements $\widetilde{g}_i \o g_i^{-1},\; 1 \leq i \leq
\widetilde{\theta},$ are central in $(U \o A)_{\sigma},$ and as in
the last part of the proof of \cite[Theorem 5.17]{AS4} it can be
seen that the map $(U \o A)_{\sigma} \to U(\D,\lambda),$
$$\;x_i \o 1 \mapsto x_i, \quad \widetilde{g}_i
\o 1 \mapsto g_i,\quad 1\o y_j \mapsto x_j,\quad 1\o g \mapsto
g$$ for all $1\leq i \leq  \widetilde{\theta}< j \leq \theta, \; g
\in \Gamma,$ induces an isomorphism of Hopf algebras
\begin{equation}\label{inductioniso}
(U \o A)_{\sigma}/(\widetilde{g}_i \o g_i^{-1} - 1 \o 1 \mid 1 \leq i \leq
\widetilde{\theta}) \cong U(\D,\lambda).
\end{equation}
Let $p_1=p_{I_1}$. By induction and Theorem \ref{R(D)}, the elements
$$x_{\beta_1}^{a_1}  \cdots x_{\beta_{p_1}}^{a_{p_1}}\widetilde{g} \o
y_{\beta_{p_1 +1}}^{a_{p_1 +1}}  \cdots
y_{\beta_{p}}^{a_{p}}g,\; a_1,\dots ,a_p \geq 0, \widetilde{g} \in
\widetilde{\Gamma}, g \in \Gamma, $$ are a basis of $U \o A.$ It
follows from \eqref{multiplicationrule} that for all $p_{1} < l \leq
p$ and $1 \leq i \leq \widetilde{\theta},$
$$(1 \o y_{\beta_l}) \cdot_{\sigma}(\widetilde{g}_i \o 1) = \chi_i(g_{\beta_l}) \widetilde{g}_i \o y_{\beta_l}.$$
Hence
$$(x_{\beta_1}^{a_1}  \cdots x_{\beta_{p_1}}^{a_{p_1}} \o
y_{\beta_{p_1 +1}}^{a_{p_1 +1}}  \cdots
y_{\beta_{p}}^{a_{p}})\cdot_{\sigma} (\widetilde{g} \o g),\; a_1,\dots a_p
\geq 0, \widetilde{g} \in \widetilde{\Gamma}, g \in \Gamma, $$ is a basis of
$(U \o A)_{\sigma}.$

Let $P = \{\widetilde{g} \o g \in (U \o A)_{\sigma} \mid \widetilde{g} \in \w{\G}, g \in
\G\}$, and let $\widetilde{P} \subset P$ be the subgroup generated
by $\widetilde{g}_i \o g_i^{-1}, \; 1 \leq i \leq \widetilde{\theta}.$ Then
$$ \G \to P/\widetilde{P},\; g \mapsto \overline{1 \o g},$$
is a group isomorphism. By \eqref{inductioniso}, $(U \o
A)_{\sigma} \o_{k[P]} k[P/\w{P}] \cong U(\D, \lambda).$ Hence
$$x_{\beta_1}^{a_1} x_{\beta_2}^{a_2} \cdots x_{\beta_p}^{a_p}g,\;
a_1, a_2,\dots, a_p \geq 0, g \in \Gamma,$$ is a basis of
$U(\D,\lambda).$

\medskip To prove (2), let $J=I_1,N = N_{J},\w\theta < i \leq
\theta$  and $\beta \in \Phi_{J}^+$. We first show that
\begin{equation}\label{rootvectorrule}
(1 \o y_i) \cdot_{\sigma} (x_{\beta}^{N} \o 1) =
\chi_{\beta}^{N}(g_i) (x_{\beta}^{N} \o 1)\cdot_{\sigma}(1 \o y_i)
\end{equation}
in $(U \o A)_{\sigma}.$ We use the notations of Section
\ref{SectionK(D)} for $\D_1$ with 
$$z^a = x_{\beta}^{N}, \text{ where }\beta = \beta_l, a = e_l \text{ for some }1 \leq l \leq p_1.$$
By \eqref{smashcomult}
$$\Delta_U(x_{\beta}^{N}) = \widetilde{g}_{\beta}^N \o x_{\beta}^{N} + x_{\beta}^N \o 1 +
\sum_{b,c \neq 0, \ub + \uc = \beta} t_{b,c}^a z^b h^c \o z^c.$$
Since $\Delta(y_i) = g_i \o y_i + y_i \o 1,$ and
$$\Delta^2(y_i) = g_i \o g_i \o y_i + g_i \o y_i \o 1 + y_i \o1 \o 1,$$
we have for all $u \in U$ by \eqref{multiplicationrule}
\begin{align}\notag
(1 \o y_i) \cdot_{\sigma} (u \o 1) &= \varphi(u\sw1)(g_i) u\sw2 \o g_i \varphi(u\sw3)(S^{-1}(y_i))\\\notag
&+ \varphi(u\sw1)(g_i) u\sw2 \o y_i \varphi(u\sw3)(1)\\\notag
&+ \varphi(u\sw1)(y_i) u\sw2 \o 1 \varphi(u\sw3)(1).\notag
\end{align}
It follows from the definition of $\varphi$ that
$$\varphi(x_{\gamma})(g) = 0 \text{ for all } \gamma \in \Phi_J^+, g \in \G.$$
Hence to compute $(1 \o y_i) \cdot_{\sigma} (u \o 1)$ with $u =
x_{\beta}^N,$ we only need to take into account the term
$\widetilde{g}_{\beta}^N \o x_{\beta}^N \o 1$ of $\Delta^2(x_{\beta}^N),$ and we
obtain
\begin{align*}
(1 \o y_i) \cdot_{\sigma} (x_{\beta}^N \o 1)&= \varphi(\widetilde{g}_{\beta}^N)(y_i\sw1) x_{\beta}^N \o y_i\sw2 \varphi(1)(S^{-1}(y_i\sw3))\\\notag
&= \varphi(\widetilde{g}_{\beta}^N)(y_i\sw1)x_{\beta}^N \o y_i\sw2\\\notag
&=\varphi(\widetilde{g}_{\beta}^N)(g_i) x_{\beta}^N \o y_i + \varphi(\widetilde{g}_{\beta}^N)(y_i)x_{\beta}^N \o 1\\\notag
&=\chi_{\beta}^{N}(g_i) (x_{\beta}^{N} \o 1)\cdot_{\sigma}(1 \o y_i),\notag
\end{align*}
since $\varphi(\widetilde{g}_{\beta}^N) | \G = \chi_{\beta}^{N}$ and
$\varphi(\widetilde{g}_{\beta}^N)(y_i) =0$ by the definition of $\varphi.$

From \eqref{inductioniso} and \eqref{rootvectorrule} we see that
for all  simple roots $\alpha \in \Phi_K^+$ with $I_1 \neq K \in \X$
and all roots $\beta \in \Phi_J^+$ with $J = I_1$
\begin{equation}\label{specialrootvector}
x_{\alpha} x_{\beta}^{N_J} = \chi_{\beta}^{N_J}(g_{\alpha})x_{\beta}^{N_J}x_{\alpha}
\end{equation}
in $U(\D,\lambda).$ Since the root vectors $x_{\alpha}$ are
homogeneous, \eqref{specialrootvector} holds for all $\alpha \in
\Phi_K^+,K \neq I_1,$ and $\beta \in \Phi_{I_1}^+.$ Since
$U(\D,\lambda)$ and the root vectors $x_{\alpha}, \alpha \in
\Phi^+,$ do not depend on the order of the connected components,
we can reorder the connected components and obtain
\eqref{specialrootvector} for all positive roots $\alpha,\beta$
lying in different connected components. For roots in the same
connected component, \eqref{specialrootvector} follows from
Theorem \ref{R(D)}. \epf

\section{ Finite-dimensional quotients}\label{Sectionfinite}

\subsection{A general criterion}\label{Sectiongeneraltheorem}
In this section we prove a generalized version of Theorem \cite[6.24]{AS5}.

Let $\G$ be an abelian group, $A$ an algebra
containing the group algebra $k[\G]$ as a subalgebra and $p \geq
1$. We assume
 $$y_1,\dots,y_p \in A,
h_1,\dots,h_p \in \G,
\psi_1, \dots, \psi_p \in \widehat{\G}, \text{ and }
N_1, \dots,N_p \geq 1,$$
such that
\begin{align}
&gy_l = \psi_l(g) y_lg, \text{ for all } 1 \leq l \leq p, g \in \G,&\label{c1}\\
&y_ky_l^{N_l} = \psi_l^{N_l}(h_k) y_l^{N_l} y_k, \text{ for all } 1 \leq k,l \leq p,&\label{c2}\\
&y_1^{a_1}\cdots y_p^{a_p}g, \; a_1,\cdots,a_p\geq 0, g \in \G, \text{ form a basis of } A.&\label{c3}
\end{align}
Let $\mathbb{T}= \{ t = (t_1,\dots,t_p) \in \Np \mid 0 \leq t_l < N_l
\text{ for all } 1 \leq l \leq p\}$.
For all $a = (a_1,\dots,a_p) \in \Np$, we define
\begin{align*}
y^a &=y_1^{a_1}\cdots y_p^{a_p},\\
\psi^a&= \psi_1^{a_1} \cdots \psi_p^{a_p},\\
aN &= (a_1N_1, \dots, a_pN_p).
\end{align*}
Hence any element $v\in A$
can be written as
\begin{equation}\label{b}
v = \sum_{t \in \mathbb{T},a \in \Np} y^t y^{aN} v_{t,a},\;
v_{t,a} \in k[\G] \text{ for all } t \in \mathbb{T},a \in \Np,
\end{equation}
where the  coefficients $v_{t,a} \in k[\G]$ are uniquely
determined.

In \cite{AS5} we assumed that $A = R \# k[\G]$, and
the subalgebra $R$ of $A$ generated by $y_1,\dots,y_p$ had the
basis $y_1^{a_1} \cdots y_p^{a_p}, a_1, \dots,a_p \geq 0.$ To see that \cite[Theorem 6.24]{AS5} extends to the more general situation considered here we first prove a more general version of \cite[Lemma 6.23]{AS5}.

We need the following commutation rules for the generators of $A$.
For all $a,b \in \mathbb{N}^p, 1 \leq l \leq p,$
\begin{align}
y^{aN}y_l&=y_ly^{aN} \psi^{aN}(h_l^{-1})\label{r1},\\
y^{bN}y^{aN} &= y^{(a+b)N}\psi^{aN}(g(b)),\label{r2}
\end{align}
where $g(b) =(g_1(b),\dots,g_p(b)) \in \G^p$ is a family of elements in $\G$ depending on $b$, and $\psi^{aN}(g(b))= \psi^{a_1N_1}(g_1(b)) \cdots \psi^{a_pN_p}(g_p(b))$. Both equations \eqref{r1} and \eqref{r2} follow from \eqref{c2}. To prove \eqref{r2} we write
\begin{align*}
y^{bN}y^{aN}&= y_1^{b_1N_1} \cdots y_p^{b_pN_p}y_1^{a_1N_1} \cdots y_p^{a_pN_p}\\
&=y_1^{b_1N_1} \cdots y_{p-1}^{b_{p-1}N_{p-1}}y_1^{a_1N_1} \cdots y_{p-1}^{a_{p-1}N_{p-1}}y_p^{(a_p+b_p)N_p}\\
&\times (\psi_1^{a_1N_1}\cdots \psi_{p-1}^{a_{p-1}N_{p-1}})(h_p^{b_pN_p})
\end{align*}
and continue in this way.

For any character $\psi \in \widehat{\G}$ let $\widetilde{\psi} : k[\G]\to k[\G]$ be the algebra map with $\widetilde{\psi}(g) = \psi(g)g$ for all $g \in \G$. Thus
$$\widehat{\G} \to \Aut(k[\G]), \; \psi \mapsto \widetilde{\psi},$$
is a group homomorphism. Then  for all $1 \leq l \leq p, a\in \mathbb{N}^p$ and $v \in k[\G]$ it follows from \eqref{c1} that
\begin{align}
vy_l &= y_l\widetilde{\psi_l}(v),\label{r3}\\
vy^{aN}&= y^{aN}\widetilde{\psi^{aN}}(v)\label{r4}.
\end{align}

\begin{Lem}\label{Lemfinite}
Let $u_l \in k[\G], 1 \leq l \leq
p$, be a family of central elements in $A$ and assume for all $1 \leq l \leq p$ that $u_l = 0$ whenever $\psi_l^{N_l} \neq \varepsilon$.
Let $M$ be a free right $k[\G]$-module with basis $m(t),t \in \mathbb{T}$,  and
define a right $k[\G]$-linear map by
$$\varphi : A \to M, \;y^t y^{aN} \mapsto m(t)u^a \text{ for all }t\in \mathbb{T},a \in \mathbb{N}^p,$$
where $u^a=u_1^{a_1} \cdots u_p^{a_p}$ for all $a = (a_1,\dots,a_p) \in \mathbb{N}^p$. Then the kernel of $\varphi$ is a right ideal of A.
\end{Lem}
\pf
Let $z \in A$ be an element of the kernel of $\varphi$. We have to show that $\varphi(zy_l)=0$ for all $1 \leq l \leq p$. We fix an arbitrary $1 \leq l \leq p$ and by \eqref{b} we have basis representations
\begin{align}
z&= \sum_{s \in \mathbb{T}, a \in \mathbb{N}^p}y^sy^{aN}v_{s,a},\label{p0}\\
y^sy_l &= \sum_{t \in \mathbb{L},b \in \mathbb{N}^p} y^t y^{bN} w_{t,b}^s \text{ for all } s \in \mathbb{T},\label{p1}
\end{align}
where the $v_{s,a}$ and the $w_{t,b}^s$ are elements in $k[\G]$.

To compute $zy_l$ we multiply \eqref{p0} with $y_l$ and then use \eqref{r3},  \eqref{r1}, \eqref{p1}, \eqref{r4} and \eqref{r2} to obtain
\begin{align*}
zy_l&= \sum_{s \in \mathbb{T}, a \in \mathbb{N}^p}y^sy^{aN}v_{s,a}y_l\\
&= \sum_{\substack{ s,t \in \mathbb{L}\\
a,b \in \mathbb{N}^p}}y^t y^{(a+b)N} \psi^{aN}(g(b))\widetilde{\psi^{aN}}(w_{t,b}^s) \psi^{aN}(h_l^{-1})\widetilde{\psi_l}(v_{s,a}).
\end{align*}
We note that $u_l^{a_l}\psi_l^{a_lN_l}(g) = u_l^{a_l}$ for all $1 \leq l \leq p, a_l \in \mathbb{N}$ and $g \in \G$, since $u_l=0$ whenever $\psi_l^{N_l} \neq \varepsilon$. Thus for all $a,b \in \mathbb{N}^p,s,t \in \mathbb{L}$
\begin{align}
u^a \psi^{aN}(g(b)) &= u^a,\; u^a \psi^{aN}(h_l^{-1}) = u^a,\label{p2}\\
u^a \widetilde{\psi^{aN}}(w_{t,b}^s) &= u^a w_{t,b}^s.\label{p3}
\end{align}
Since $u^a$ is central in $A$ it follows from \eqref{c3} that
\begin{equation}\label{p4}
u^a \widetilde{\psi_l}(v_{s,a}) = \widetilde{\psi_l}(u^av_{s,a}).
\end{equation}
Hence
\begin{align*}
\varphi(zy_l)&= \sum_{\substack{ s,t \in \mathbb{T}\\
a,b \in \mathbb{N}^p}} m(t) u^{a+b} \psi^{aN}(g(b))\widetilde{\psi^{aN}}(w_{t,b}^s) \psi^{aN}(h_l^{-1})\widetilde{\psi_l}(v_{s,a})\\
&=\sum_{\substack{ s,t \in \mathbb{T}\\
b \in \mathbb{N}^p}} m(t)u^b w_{t,b}^s \widetilde{\psi_l}(\sum_{a \in \mathbb{N}^p} u^a v_{s,a}),
\end{align*}
where the second equality follows from \eqref{p2}, \eqref{p3} and \eqref{p4}. This proves our claim since $\varphi(z)=0$, and therefore
$$\sum_{a \in \mathbb{N}^p} u^a v_{s,a}=0 \text{ for all } s \in \mathbb{T}.$$
\epf

\begin{Thm}\label{Theoremfinite}
Assume the situation above. Let $u_l \in k[\G], 1 \leq l \leq
p$, be a family of elements in the group algebra. Then the following are equivalent:
\begin{enumerate}
\item The residue classes of $y^tg, \; t \in \mathbb{T}, g \in \G,$ form
a basis of the quotient algebra $A/(y_l^{N_l} - u_l \mid 1 \leq l
\leq p).$
\item For all $1 \leq l \leq p$, $u_l$ is central in $A$, and if
$\psi_l^{N_l} \neq \varepsilon$, then $u_l = 0.$
\end{enumerate}
\end{Thm}
\pf As in the proof of \cite[Theorem 6.24]{AS5} this follows from Lemma \ref{Lemfinite}.
\epf

\subsection{The Hopf algebra $u(\D,\lambda,\mu)$}\label{Sectionu}

Let $\G$ be a finite abelian group, and $\D = \D(\Gamma, (g_i)_{1
\leq i \leq \theta}, (\chi_i)_{1 \leq i \leq \theta},
(a_{ij})_{1\leq i,j \leq \theta})$  a datum of finite Cartan type.
We assume the situation of Section \ref{Sectionnotations}.

\begin{Def}\label{Defu(D)}
A family $\mu=(\mu_{\alpha})_{\alpha \in \Phi^+}$ of elements in
$k$ is called a {\em family of root vector parameters for $\D$} if
the following condition is satisfied for all $\alpha \in
\Phi_J^+,J \in \X$:  If $g_{\alpha}^{N_J} =1$ or
$\chi_{\alpha}^{N_J} \neq \varepsilon,$ then $\mu_{\alpha} =0.$

Let $\mu$ be a family of root vector parameters for $\D$. For all
$J \in \X$, and $\alpha \in \Phi_J^+,$ we define
\begin{equation}\label{ualphageneral}
\pi_J(\mu) = (\mu_{\beta})_{\beta \in \Phi_J^+}, \text{ and }
u_{\alpha}(\mu) = u_{\alpha}(\pi_J(\mu)),
\end{equation}
where $u_{\alpha}(\pi_J(\mu))$ is introduced in Definition \ref{ualpha}.
Let $\lambda$ be a family of linking parameters for $\D$. Then we define
\begin{equation}\label{Definitionu(D)}
u(\D,\lambda,\mu)=U(\D,\lambda)/(x_{\alpha}^{N_J} -
u_{\alpha}(\mu) \mid \alpha \in \Phi_J^+, J \in \X).
\end{equation}
\end{Def}
By abuse of language we still write $x_i$ and $g$ for the images
of $x_i$ and $g \in \G$ in $u(\D,\lambda,\mu).$ For all $1 \leq l
\leq p$, we define $N_l = N_J$, if $\beta_l \in \Phi_J^+, J \in
\X.$

\begin{Lem}\label{Lemcentral}
Let $\D, \lambda$ and $ \mu$ as above, and $\alpha \in \Phi^+.$
Then $u_{\alpha}(\mu)$ is central in $U(\D,\lambda).$
\end{Lem}
\pf Let $\alpha \in \Phi_J^+$, where $J \in \X,$ and $N = N_J.$
To simplify the notation, we assume $J = I_1 =
\{1,2,\dots,\w{\theta}\},$ and $\Phi_J^+ =
\{\beta_1,\beta_2,\dots,\beta_{p_1}\}.$ We apply the results and
notations of Section \ref{SectionK(D)} to the connected component
$I_1$. Let $\widetilde{\mu} = (\mu_l)_{1 \leq l \leq p_1}, \D_1 = \D(\G, (g_i)_{1 \leq i \leq \w{\theta}},(\chi_i)_{1 \leq i \leq \w{\theta}},(a_{i,j})_{1 \leq i,j \leq \w{\theta}})$, and 
$$\varphi_{\widetilde{\mu}} : K(\D_1) \# k[\G] \to k[\G]$$
the Hopf algebra map defined by $\widetilde{\mu}$ in Theorem \ref{mainconstruction}.
As in Section \ref{SectionK(D)} we define $u^a = \varphi_{\widetilde{\mu}}(z^a)$ for all $0 \neq a \in \mathbb{N}^{p_1}$. Thus $u_{\alpha}(\mu)= u^{e_l}$ where $\alpha = \beta_l$.

We show by induction on $\het(\ua)$ that 
\begin{equation}\label{centralua}
x_i u^a = u^a x_i \text{ for all } 1 \leq i \leq \theta, 0 \neq a \in \mathbb{N}^{p_1}.
\end{equation}
To prove \eqref{centralua} we can assume that $u^a \neq 0$. By \eqref{ua} it suffices to show that
\begin{equation}\label{centralmu}
x_i h^a = h^a x_i \text{ for all } 1 \leq i \leq \theta, 0 \neq a \in \mathbb{N}^{p_1} \text{ with }u^a \neq 0.
\end{equation}
Recall that $h^a = g_{\beta_1}^{Na_1} \cdots g_{\beta_{p_1}}^{Na_{p_1}}$ for $a = (a_1,\dots,a_{p_1})$.

Let $1 \leq i \leq \theta$ and $0 \neq a \in \mathbb{N}^{p_1} \text{ with }u^a \neq 0$.
 Let $1 \leq l \leq
p_1$ and $\beta_l = \sum_{j=1}^{\w{\theta}} n_j \alpha_j,$
where $n_j \in \mathbb{N}$ for all $1 \leq j \leq \w{\theta}.$
Then by definition, $g_{\beta_l} = \prod_{1 \leq j \leq
\w{\theta}} g_j^{n_j},$ and $\chi_{\beta_l} = \prod_{1 \leq j \leq
\w{\theta}} \chi_j^{n_j}.$ Hence
$$\chi_i(g_{\beta_l}^N) \chi_{\beta_l}^N(g_i) = \prod_{1 \leq j \leq
\w{\theta}} q_{ii}^{a_{ij}Nn_j} =1,$$ since $q_{ii}^N = 1$, if $i
\in I_1,$ and $a_{ij} =0$, if $i \notin I_1.$ Since $u^a \neq 0$, it follows from Lemma \ref{partial} that  $\chi_{\beta_l}^N = \varepsilon$ for all $1 \leq l
\leq p_1$ with $a_l >0.$ Hence $\chi_i(g_{\beta_l}^N) =1$
for all $l$ with $a_l>0.$ This implies \eqref{centralmu} since 
$$h^a
x_i = \chi_i(h^a) x_i h^a.$$
\epf

\begin{Thm}\label{Theoremu(D)}
Let $\D$ be a datum of finite Cartan type satisfying
\eqref{orderqij} and \eqref{orderNJ}. Let $\lambda$ and $\mu$ be
families of linking and root vector parameters for $\D$. Then
$u(\D,\lambda,\mu)$ is a quotient Hopf algebra of $U(\D,\lambda)$
with group-like elements $G(u(\D,\lambda,\mu)) \cong \G,$ and the
elements
$$x_{\beta_1}^{a_1}x_{\beta_2}^{a_2} \cdots x_{\beta_p}^{a_p}g,\;  0
\leq a_l < N_l,\; 1 \leq l \leq p,\; g \in \G$$ form a  basis of
$u(\D,\lambda,\mu).$ In particular,
$$\dim u(\D,\lambda,\mu)= \prod_{J \in  \X} N_J^{|\Phi_J^+|} |\G|.$$
\end{Thm}
\pf
By Theorem \ref{TheoremU(D)}, the elements
$$x_{\beta_1}^{a_1}x_{\beta_2}^{a_2} \cdots x_{\beta_p}^{a_p}g,\;  0
\leq a_l,\;1 \leq l \leq p,\; g \in \G$$ are a basis of
$U(\D,\lambda)$. We want to apply Theorem \ref{Theoremfinite} to $U(\D,\lambda)$ and 
$$y_l = x_{\beta_l},\; h_l=g_{\beta_l},\; \psi_l = \chi_{\beta_l},\;u_l =
u_{\beta_l}(\mu),\;1 \leq l \leq p.$$ 
Conditions \eqref{c2}, \eqref{c3} are satisfied by Theorem \ref{TheoremU(D)}. To check the conditions in Theorem \ref{Theoremfinite} (2) we apply for each connected component $J \in \X$
the results of Section \ref{SectionK(D)} with
$$\eta_l = \chi_{\beta_l}^{N_l},\; 1 \leq l \leq p,\; \beta_l \in \Phi_J^+.$$
Since $\varphi_{\mu}$ is a Hopf algebra map by Theorem \ref{mainconstruction} it follows from Lemma \ref{partial} that $u_{\beta_l}(\mu) = 0$ if $\chi_{\beta_l}^{N_l} \neq \varepsilon$. By Lemma
\ref{Lemcentral}, $u_{\beta_l}(\mu)$ is central in $U(\D,\lambda).$
Hence  the claim about the basis of $u(\D,\lambda,\mu)$
follows from Theorem
\ref{Theoremfinite}.

\medskip
We now show that $u(\D,\lambda,\mu)$ is a Hopf algebra. Let $J \in
\X$. We denote the restriction of $\D$ to the connected component
$J$ by $\D_J$. By Theorem \ref{mainconstruction}, the map
$\varphi_{\mu} : K(\D_J) \# k[\G] \to k[\G]$ is a Hopf algebra
homomorphism. The kernel of $\varphi_{\mu}$ is generated by all
$x_{\alpha}^{N_J} - u_{\alpha}(\mu)$ with  $\alpha \in \Phi_J^+.$ Hence
the elements $x_{\alpha}^{N_J} - u_{\alpha}(\mu), \alpha \in
\Phi_J^+,$ generate a Hopf ideal in $K(\D_J) \# k[\G]$ and in
$U(\D,\lambda).$

\medskip

The Hopf algebra $u(\D,\lambda,\mu)$ is generated by the
skew-primitive elements $x_1,\dots,x_\theta$ and the image of
$\G$. Hence $G(u(\D,\lambda,\mu)) \cong \G.$ \epf

For explicit examples of the Hopf algebras $u(\D,\lambda,\mu)$ see
\cite[Section 6]{AS5} for type $A_n,n \geq 1,$ and \cite{BDR} for
type $B_2.$ In these papers, and for these types, the elements
$u_{\alpha}(\mu)$ are precisely written down. It is an interesting
problem to find an explicit algorithm describing the
elements $u_{\alpha}(\mu)$ for any connected Dynkin diagram.

\section{The associated graded Hopf algebra}\label{Sectiongraded}

\subsection{Nichols algebras}\label{SectionNichols}

To determine the structure of a given pointed Hopf algebra, we
proceed as in \cite{AS1} and study the associated graded Hopf
algebra.

Let $A$ be a pointed Hopf algebra with group of group-like
elements $G(A) = \G.$ Let
$$A_0 =k[\G] \subset A_1 \subset \dots \subset A,\; A=\cup_{n \geq 0}A_n$$
be the coradical filtration of $A$. We define the associated
graded Hopf algebra \cite[5.2.8]{M} by
$$\gr(A) = \oplus_{n\geq0} A_n/A_{n-1},\; A_{-1}=0 .$$

Then $\gr(A)$ is a pointed Hopf algebra with the same dimension
and coradical as $A$. The projection map $\pi : \gr(A) \to k[\G]$
and the inclusion $\iota : k[\G] \to \gr(A)$ are Hopf algebra maps
with $\pi\iota  = \id_{k[\G]}.$ Let
\begin{equation}\label{DefR}
R = \{ x \in \gr(A) \mid (\id \o \pi)\Delta(x) = x \o 1\}
\end{equation}
be the algebra of $k[\G]$-coinvariant elements. Then
$R=\oplus_{n\geq0} R(n)$ is a graded Hopf algebra in $\YDG$, and
by \eqref{Radfordiso}
\begin{equation}\label{writeR}
\gr(A) \cong R \# k[\G].
\end{equation}
Let $V=P(R) \in {\YDG}$ be the Yetter-Drinfeld module of primitive
elements in $R.$ We call its braiding
$$c : V \o V \to V \o V$$
the {\em infinitesimal braiding of $A$}.

Let $\mathfrak{B}(V)$ be the subalgebra of $R$ generated by $V.$
Thus $B= \mathfrak{B}(V)$ is the {\em Nichols algebra} of $V$
\cite{AS2}, that is,
\begin{align}
&B= \oplus_{n\geq0} B(n) \text{ is a graded Hopf algebra in } \YDG,&\\
&B(0) = k 1,\;  B(1) =V,&\\
&B(1) = P(B),&\\
&B \text{ is generated as an algebra by B(1)}.&
\end{align}
$\mathfrak{B}(V)$ only depends on the vector space $V$  with its
Yetter-Drinfeld structure (see the discussion in \cite[Section
2]{AS5}). As an algebra and coalgebra, $\mathfrak{B}(V)$ only
depends on the braided vector space $(V,c).$

We assume in addition that $A$ is finite-dimensional and $\G$ is
abelian. Then there are  $g_1,\dots,g_{\theta} \in \G,$
$\chi_1,\dots,\chi_{\theta} \in \widehat{\G}$ and a basis
$x_1,\dots,x_{\theta}$ of $V$ such that $x_i \in V_{g_i}^{\chi_i}$
for all $ 1 \leq i \leq \theta.$ We call
$$(q_{ij}=\chi_j(g_i))_{1 \leq i,j \leq \theta}$$
the {\em infinitesimal braiding matrix of $A$}.

A braiding matrix $(q_{ij})_{1 \leq i,j \leq \theta}$ whose entries $q_{ij}$ are roots of unity is of {\em Cartan type} if $q_{ii} \neq 1$ for all $1 \leq i \leq \theta$, and if there are integers $a_{ij}, 1 \leq i,j \leq \theta$, such that for all $1 \leq i,j \leq \theta$
$$q_{ij}q_{ji}=q_{ii}^{a_{ij}}.$$
We can assume that $a_{ii}=2$ for all $1 \leq i \leq \theta$, and
$$-\ord(q_{ii}) < a_{ij} \leq 0 \text{ for all } 1 \leq i,j \leq \theta.$$ Then the matrix $(a_{ij})$ is uniquely determined. It is a generalized Cartan matrix and is called the {\em Cartan matrix of $(q_{ij})$} \cite{AS2}.

The first step to classify pointed Hopf algebras is the
computation of the Nichols algebra. We begin with the description of Nichols algebras of Yetter-Drinfeld modules of finite Cartan type.

\begin{Thm}\label{CartanNichols}
Let $\D = \D(\Gamma, (g_i)_{1 \leq i \leq \theta}, (\chi_i)_{1
\leq i \leq \theta}, (a_{ij})_{1\leq i,j \leq \theta})$  be a
datum of finite Cartan type with finite abelian group $\G$. Assume
\eqref{orderqij} and \eqref{orderNJ}. Let $V \in {\YDG}$ be a
vector space with basis $x_1,\dots,x_{\theta}$ and $x_i \in
V_{g_i}^{\chi_i}$ for all $1 \leq i \leq \theta.$ Then
$\mathfrak{B}(V)$ is the quotient algebra of $T(V)$ modulo the
ideal generated by the elements
\begin{align}
&\ad_c(x_i)^{1 - a_{ij}}(x_j) \text{ for all } 1 \leq i,j \leq \theta, i \neq j,&\label{SerreNichols}\\
&x_{\alpha}^{N_J} \text{ for all } \alpha \in \Phi_J^+, J \in \mathcal{X}.&\label{rootNichols}
\end{align}
\end{Thm}
\pf
Using results of Lusztig \cite{L1},\cite{L2}, Rosso \cite{Ro} and
M\"uller \cite{M1} and twisting we proved this theorem in \cite[Theorem
4.5]{AS4} assuming in addition that $\ord(q_{ij})$ is odd for all $1 \leq i,j \leq \theta, i \neq j$. By Lemma \ref{q} the proof of \cite[Theorem
4.5]{AS4} works without this additional assumption.
\epf

\begin{Cor}\label{gru}
Assume the situation of Theorem \ref{CartanNichols}, and let
$\lambda$ and $\mu$ be linking and root vector parameters for
$\D$. Then
$$\gr(u(\D,\lambda,\mu)) \cong u(\D,0,0) \cong \mathfrak{B}(V) \# k[\G].$$
\end{Cor}
\pf
Let $A=u(\D,\lambda,\mu).$ There is a well-defined Hopf algebra map
$$u(\D,0,0) \to \gr(u(\D,\lambda,\mu)),$$
mapping $x_i,1 \leq i \leq \theta,$ onto the residue class of
$x_i$ in $A_1/A_0,$ and $g \in \G$ onto $g.$ Since
$\dim(u(\D,0,0)) =
\dim(u(\D,\lambda,\mu))=\dim(\gr(u(\D,\lambda,\mu))$ by Theorem
\ref{Theoremu(D)}, it follows that $u(\D,0,0) \cong
\gr(u(\D,\lambda,\mu))$. By Theorem \ref{CartanNichols},
$u(\D,0,0)  \cong \mathfrak{B}(V)\# k[\G].$ \epf

\bigbreak In \cite{AS2} and \cite{AS4} we determined the structure
of finite-dimensional Nichols algebras assuming that $V$ is of
Cartan type and satisfies some more assumptions in the case of
small orders ($\leq 17$) of the diagonal elements $q_{ii}$. Recent
results of Heckenberger \cite{H1}, \cite{H2}, \cite{H3} together
with Theorem \ref{CartanNichols} allow to prove the following very
general structure theorem on Nichols algebras.

\begin{Thm}\label{classificationNichols}
Let $\G$ be a finite abelian group, and $V \in {\YDG}$ a
Yetter-Drinfeld module such that $\mathfrak{B}(V)$ is
finite-dimensional. Choose a basis $x_i \in V$ with $x_i \in
V_{g_i}^{\chi_i}, g_i \in \G,\chi_i \in \widehat{\G},  \text{ for
all } 1 \leq i \leq \theta.$ For all $1 \leq i,j \leq \theta,$
define $q_{ij}= \chi_j(g_i),$ and assume
\begin{align}
&\ord(q_{ii}) \text{ is odd},\label{H1}\\
&\ord(q_{ii}) \text{ is prime to 3 if $q_{il}q_{li}\in \{q_{ii}^{-3},q_{ll}^{-3}\}$ for some $l$},\label{G2}\\
&\ord(q_{ii})>3 .\label{H2}
\end{align}
Then there is a datum $\D = \D(\G,(g_i)_{1\leq i \leq
\theta},(\chi_i)_{1\leq i \leq \theta},(a_{ij})_{1\leq i,j \leq
\theta})$ of finite Cartan type such that
$$ \mathfrak{B}(V) \# k[\G] \cong u(\D,0,0).$$
\end{Thm}
\pf Since $\mathfrak{B}(V)\# k[\G]$ is finite-dimensional, $q_{ii}
\neq 1$ for all $1 \leq i \leq \theta$ by \cite[Lemma 3.1]{AS1}.

For all $1 \leq i,j \leq \theta, i \neq j,$ let $V_{ij}$ be
the vector subspace of $V$ spanned by $x_i,x_j.$ Then
$\mathfrak{B}(V_{ij})$ is isomorphic to a subalgebra of
$\mathfrak{B}(V),$ hence it is finite-dimensional. Heckenberger
\cite{H1}, \cite{H2} classified finite-dimensional Nichols
algebras of rank 2. By \eqref{H1} and \eqref{H2} it follows from the list in
\cite[Theorem 4]{H1} that $V_{ij}$ is of finite Cartan type, that
is, there are $a_{ij}, a_{ji} \in \{0,-1,-2,-3\}$ with $a_{ij}
a_{ji} \in \{0,1,2,3\},$ and
$$q_{ij}q_{ji}=q_{ii}^{a_{ij}}=q_{jj}^{a_{ji}}.$$
Thus $(q_{ij})_{1 \leq i,j \leq \theta}$ is of Cartan type in the
 with generalized Cartan matrix
$(a_{ij})$. In \cite[Theorem 4]{H3} Heckenberger extended part
(ii) of \cite[Theorem 1.1]{AS2} (where we had to exclude some
small primes) and showed that a diagonal braiding $(q_{ij})$ of a
braided vector space $V$ is of finite Cartan type if it is of
Cartan type and $\mathfrak{B}(V)$ is finite-dimensional. Hence
$(a_{ij})$ is a Cartan matrix of finite type, and the claim
follows from Theorem \ref{CartanNichols}. \epf

\subsection{Generation in degree one}\label{SectionGeneration}
We generalize our results in \cite[Section 7]{AS4}. Let $A$ be a
finite-dimensional pointed Hopf algebra with $\G,V,$ and $R$ as in
Section \ref{SectionNichols}. To prove that $\mathfrak{B}(V)=R$,
we dualize. Let $S =R^*$ the dual Hopf algebra in $\YDG$ as in
\cite[Lemma 5.5]{AS2}. Then $S = \oplus_{n \geq0} S(n)$ is a
graded Hopf algebra in $\YDG$, and by \cite[Lemma 5.5]{AS2}, $R$
is generated in degree one, that is, $\mathfrak{B}(V) = R$, if and
only if $P(S) = S(1).$ The dual vector space $S(1)$ of $V=R(1)$
has the same braiding $(q_{ij})$ (with respect to the dual basis)
as $V$. Our strategy to show $P(S) = S(1)$ is to identify $S$ as a
Nichols algebra.
In the next Lemma we use \cite{H1, H2} to prove a very general
version of \cite[Lemma 7.2]{AS4}.

\begin{Lem}\label{GenerationSerre}
Let $\D = \D(\G,(g_i)_{1\leq i \leq \theta}, (\chi_i)_{1\leq i
\leq \theta},(a_{ij})_{1\leq i,j \leq \theta})$ be a datum of
finite Cartan type with finite abelian group $\G$. Let $S =
\oplus_{n\geq0} S_n$ be a finite-dimensional graded Hopf algebra
in $\YDG$ with $S(0) = k1,$ and let $x_1,\dots,x_{\theta}$ be a
basis of $S(1)$ with $x_i \in S(1)_{g_i}^{\chi_i}$ for all $1 \leq
i \leq \theta.$ Assume \eqref{H1} and
\begin{equation}\label{H3}
\ord(q_{ii}) >7 \text{ for all }1 \leq i \leq \theta.
\end{equation}
Then
\begin{align}
&\ad_c(x_i)^{1 - a_{ij}}(x_j)=0 \text{ for all } 1 \leq i,j \leq
\theta,\, i \neq j.&\label{SSerre}
\end{align}
\end{Lem}
\pf We first note that the Nichols algebra of the primitive
elements $P(S) \in {\YDG}$ is finite-dimensional. This can be seen
by looking at $\gr(S \# k[\G]).$

Assume that there are $1 \leq i,j \leq \theta, i \neq j,$ with
$\ad_c(x_i)^{1 - a_{ij}}(x_j)\neq 0.$ We define
$$y_1 = x_i,\quad y_2=\ad_c(x_i)^{1 - a_{ij}}(x_j).$$
By \cite[A.1]{AS2}, $y_2$ is a primitive element. Since $y_1,y_2$
are non-zero elements of different degree, they are linearly
independent. We know that the Nichols algebra of $W=ky_1 +ky_2$ is
finite-dimensional, since $\mathfrak{B}(P(S))$ is
finite-dimensional. We denote
$$h_1 = g_i,\, h_2 = g_i^{1-a_{ij}}g_j \in \G, \text{ and } \eta_1 =
\chi_i, \,\eta_2 = \chi_i^{1 -a_{ij}} \chi_j \in \widehat{\G}.$$
Thus $y_i \in S_{h_i}^{\eta_i}, 1 \leq i \leq 2.$ Let
$(Q_{ij}=\eta_j(h_i))_{1 \leq i,j \leq 2}$ be the braiding matrix
of $y_1$, $y_2.$ We compute
$$Q_{11} = q_{ii},\quad Q_{22} = q_{ii}^{1- a_{ij}} q_{jj},\quad
Q_{12}Q_{21} = q_{ii}^{2 - a_{ij}}.$$ By assumption, the order of
$Q_{11}=q_{ii}$ is odd and $>3$. Since $\mathfrak{B}(W)$ is
finite-dimensional, $Q_{22} \neq 1$ by \cite[Lemma 3.1]{AS1}. Thus
$Q_{22}$ has odd order, since the orders of $q_{ii},q_{jj}$ are
odd. By checking Heckenberger's list in \cite[Theorem 4]{H1}, and
thanks to \cite{H2}, we see that the braiding $(Q_{ij})$ is of
finite Cartan type or that we are in case (T3) with
$$Q_{12}Q_{21} = Q_{11}^{-1}.$$
Hence there exists $A_{12} \in \{0,-1,-2,-3\}$ with
$$Q_{12}Q_{21} = Q_{11}^{A_{12}}.$$
Since $Q_{12}Q_{21}=q_{ii}^{2 - a_{ij}},$ and $Q_{11} = q_{ii},$
it follows that the order of $q_{ii}$ divides $2-a_{ij}-A_{12} \in
\{2,3,4,5,6,7,8\}.$ This is a contradiction since the order of
$q_{ii}$ is odd and $>7.$ \epf


The next theorem is one of the main results of this paper.

\begin{Thm}\label{TheoremGeneration}
Let $A$ be a finite-dimensional pointed Hopf algebra with abelian
group $G(A) = \G$ and infinitesimal braiding matrix $(q_{ij})_{1
\leq i,j \leq \theta}.$ Assume \eqref{H1}, \eqref{G2} and \eqref{H3}.
Then $A$ is generated by group-like and skew-primitive elements,
that is,
$$R = \mathfrak{B}(V),$$
where $R$ is defined by \eqref{DefR}, and $V=R(1).$
\end{Thm}
\pf We argue as in the proof of \cite[Theorem 7.6]{AS4}. Let $S =
R^*$ be the dual Hopf algebra in $\YDG.$ Then $S(1)=R(1)^*$ has
the same braiding $(q_{ij})$ as $R(1)$ with respect to the dual
basis $(x_i)$ of the corresponding basis of $R(1).$ By Theorem
\ref{classificationNichols} $(q_{ij})$ is of finite Cartan type.
By Lemma \ref{GenerationSerre} the Serre relations
\eqref{SerreNichols} hold for the elements $x_i.$ Then the root
vector relations \eqref{rootNichols} follow by \cite[Lemma
7.5]{AS4}. Hence $S \cong \mathfrak{B}(S(1))$ by Theorem
\ref{CartanNichols}, and $S(1) = P(S).$ By duality, $R$ is a
Nichols algebra.
\epf

\section{Lifting}\label{SectionLifting}

From Section \ref{Sectiongraded} we know a presentation of
$\gr(A)$ by generators and relations under the assumptions of
Theorems \ref{classificationNichols} and \ref{TheoremGeneration}.
To lift this presentation to $A$ we need the following formulation
of \cite[Lemma 5.4]{AS1} which is a consequence of the theorem of
Taft and Wilson \cite[Theorem 5.4.1]{M}. Here it is crucial that
the group is abelian.

\begin{Lem}\label{LemmaLifting}
Let $A$ be a finite-dimensional pointed Hopf algebra with abelian
group $G(A) = \G.$ Write $\gr(A) \cong R \# k[\G]$ as in
\eqref{writeR}, and let $V = R(1)$ with basis $x_i \in
V_{g_i}^{\chi_i},  g_i \in \G,\chi_i \in \widehat{\G},  1 \leq i
\leq \theta.$ Let $A_0 \subset A_1$ be the first two terms of the
coradical filtration of $A$. Then
\begin{align}
&\oplus_{g,h \in \G, \varepsilon \neq \chi \in \widehat{\G}} P_{g,h}^{\chi}(A) \xrightarrow{\cong} A_1/A_0 \xleftarrow{\cong} V \# k[\G].&\label{liftingxi}\\
& \text{For all } g \in \G,P_{g,1}(A)^{\varepsilon} = k (1-g), \text{ and if } \varepsilon \neq \chi \in \widehat{\G},\text{ then }&\label{skewprimitive1}\\
&P_{g,1}(A)^{\chi} \neq 0 \iff g=g_i, \chi=\chi_i, \text{ for some } 1 \leq i \leq \theta.&\label{skewprimitive2}
\end{align}
\end{Lem}

We can now prove our main structure theorem.

\begin{Thm}\label{Theorempointed}
Let $A$ be a finite-dimensional pointed Hopf algebra with abelian
group $G(A) = \G$ and infinitesimal braiding matrix $(q_{ij})_{1
\leq i,j \leq \theta}.$ Assume \eqref{H1}, \eqref{G2} and
\eqref{H3}. Then
$$A \cong u(\D,\lambda,\mu),$$
where $\D =\D(\Gamma, (g_i)_{1 \leq i \leq \theta}, (\chi_i)_{1
\leq i \leq \theta}, (a_{ij})_{1\leq i,j \leq \theta})$ is a datum
of finite Cartan type, and $\lambda$ and $\mu$ are families of
linking and root vector parameters for $\D.$
\end{Thm}
\pf By Theorems \ref{classificationNichols} and
\ref{TheoremGeneration}, there is a datum $\D$ of finite Cartan
type such that $\gr(A) \cong u(\D,0,0).$ By Lemma
\ref{LemmaLifting}, for all $1 \leq i \leq \theta$ we can choose
$$a_i \in P(A)_{g_i,1}^{\chi_i} \text{ corresponding to } x_i \text{ in \eqref{liftingxi}}.$$
We have shown in Theorem \cite[6.8]{AS4} that
\begin{align*}
&\ad_c(a_i)^{1-a_{ij}}(a_j)=0, \text{ for all } 1\leq i,j \leq \theta, i \sim j, i \neq j,&\\
&a_i a_j - q_{ij} a_j a_i - \lambda_{ij}(1 - g_i g_j)=0, \text{ for all } 1 \le i < j \le \theta, \, i\not\sim j,&
\end{align*}
for some family $\lambda$ of linking parameters. Thus there is a
homomorphism of Hopf algebras
$$\varphi : U(\D,\lambda) \to A,\;\varphi| \G = \id_{\G},\,
\varphi(x_i) =  a_i, \text{ for all } 1 \leq i \leq \theta.$$ By
Theorem \ref{TheoremGeneration}, $\varphi$ is surjective.

\medskip

We now use the notation of Section \ref{SectionK(D)} and show that
\begin{align}\label{Lemmavarphi}
&\varphi(x_{\alpha}^{N_J}) \in k[\G] \text{ for all } \alpha \in \Phi_J^+, J \in \mathcal{X}.&
\end{align}
We fix $J \in \mathcal{X}$ with $p = |\Phi_J^+|,$ and show by induction on $\het(\ua)$ that
\begin{align}\label{Lemmavarphigen}
&\varphi(z^a) \in k[\G] \text{ for all } a \in \Np.&
\end{align}
Let $0\neq a \in \Np.$ Since $\varphi$ is a Hopf algebra map, we
see from \eqref{smashcomult} that
$$\Delta(\varphi(z^a)) = h^a \o \varphi(z^a) + \varphi(z^a) \o 1 + w,$$
where by induction
$$w = \sum_{b,c \neq o, \ub + \uc = \ua} t_{b,c}^a \,\varphi(z^b)
h^c \o \varphi(z^c) \in k[\G] \o k[\G].$$ In particular,
$\varphi(z^a) \in A_1$ by definition of the coradical filtration.
We multiply this equation with $g\o g, g \in \G,$ from the left
and $g^{-1} \o g^{-1}$ from the right. Since $gz^ag^{-1} =
\eta^a(g) z^a$, we obtain $w = \eta^a(g) w$ for all $g \in \G.$

Suppose $\eta^a \neq \varepsilon.$ Then $w =0,$ and $\varphi(z^a)
\in P_{h^a,1}^{\eta^a}.$ Then $\varphi(z^a) =0$ by Lemma
\ref{LemmaLifting} \eqref{skewprimitive2}, since $\chi_l(g_l) \neq
1$ for all $1 \leq l \leq \theta,$ but $\eta^a(h^a) = 1$ by the
Cartan condition (see the proof of \cite[Lemma 7.5]{AS2} for a
similar computation).

If $\eta^a =\varepsilon,$ then $\varphi(z^a) \in A_1^{\varepsilon}
= k[\G]$ by Lemma \ref{LemmaLifting} \eqref{skewprimitive1}.

\medskip
This proves \eqref{Lemmavarphigen} and \eqref{Lemmavarphi}. Then
we conclude for each $J \in \mathcal{X}$ from Theorem
\ref{mainconstruction} that the map
$$K(\D_J) \# k[\G] \to U(\D, \lambda) \xrightarrow{\varphi} A$$
has the form $\varphi_{\mu^J}$ for some family of scalars $\mu^J$
as in Theorem \ref{mainconstruction} for the connected component
$J.$ Define $\mu = (\mu_{\alpha})_{\alpha \in \Phi+}$ by
$\mu_{\alpha} = \mu_{\alpha}^J$ for all $ \alpha \in \Phi_J^+.$
Then $\mu$ is a family of root vector parameters for $\D,$ and the
elements $u_{\alpha}(\mu) \in k[\G]$ are defined in
\eqref{ualphageneral} for each $J \in \mathcal{X}$ and $\alpha \in
\Phi_J^+.$ It follows that $\varphi(x_{\alpha}^{N_J}) =
u_{\alpha}(\mu) = \varphi(u_{\alpha}(\mu))$ for all $J \in
\mathcal{X}, \alpha \in \Phi_J^+.$ Thus $\varphi$ factorizes over
$u(\D,\lambda,\mu)$. Since 
$$\dim(A) = \dim(\gr(A)) =
\dim(u(\D,0,0)) = \dim(u(\D,\lambda,\mu))$$ by Theorem
\ref{Theoremu(D)}, $\varphi$ induces an isomorphism
$u(\D,\lambda,\mu) \cong A.$ \epf
\begin{Cor}\label{Cauchy}
Let $A$ be a finite-dimensional pointed Hopf algebra with abelian
group $G(A) = \G$ satisfying the assumptions of Theorem
\ref{Theorempointed}. Then for each prime divisor $p$ of the
dimension of $A$ there is a group-like element of order $p$ in
$A$.
\end{Cor}
\pf This follows from Theorems \ref{Theorempointed} and
\ref{Theoremu(D)}. \epf We note that the analog of Cauchy's
theorem in group theory is false for arbitrary, non-pointed Hopf
algebras. Let $A$ be a finite-dimensional Hopf algebra with only
trivial group-like elements, such as the dual of the group algebra
of a finite group $G$ with $G = [G,G].$ Then $A$ does not contain
any Hopf subalgebra of prime dimension, since any Hopf algebra of
prime dimension is a group algebra by Zhu's theorem \cite{Z}.

Cauchy's theorem for semisimple Hopf algebras in a version conjectured by Etingof and Gelaki was recently shown in \cite{KSZ}: Each prime divisor of a semisimple Hopf algebra divides the exponent of the Hopf algebra.

\section{Isomorphism classes}\label{SectionIso}
In this last section we determine all isomorphisms between the
Hopf algebras $u(\D,\lambda,\mu)$ in terms of some universal constants. We explicitly computed these constants for connected components of type $A$ in \cite{AS7}.

For convenience we introduce a normalization condition for Cartan
matrices and their root systems. Let $(a'_{ij})_{1 \leq i,j \leq
\theta}$ and $(a_{ij})_{1 \leq i,j \leq \theta}$ be Cartan
matrices of finite type. A {\em diagram isomorphism} between
$(a'_{ij})$ and $(a_{ij})$ is a permutation $\tau$ of
$\{1,2,\dots,\theta\}$ with $a'_{ij}= a_{\tau(i),\tau(j)}$ for all
$1 \leq i,j \leq \theta.$ We choose from each isomorphism class of
connected Cartan matrices of finite type one representative. The
chosen representatives are called {\em standard} Cartan matrices.
We fix a reduced representation of the longest element in the Weyl
group and the corresponding ordering of the positive roots of the
standard matrices as described in Section \ref{Cartan}.

From now on we assume that any Cartan datum $\D $ satisfies the
following {\em additional normalizing condition}:
\begin{enumerate}
\item[]
The Cartan matrix $(a_{ij})_{1 \leq i,j \leq \theta}$ of $\D$ is a
block diagonal matrix, and each matrix on the diagonal is one of
the standard connected Cartan matrices. Moreover for each
connected component $J$ of $I$ we fix the same order of the
positive roots as for the chosen representative.
\end{enumerate}
Thus up to a shift of indices we can identify the Cartan matrix of
any connected component with a standard Cartan matrix. We use the
fixed ordering of the positive roots in each component to define
the root vectors of $\D$.

Note that any diagram isomorphism induces an isomorphism of the
corresponding Nichols algebras. Hence up to Hopf algebra
isomorphisms we can assume by the proof of Theorem
\ref{Theorempointed} that the normalizing condition is satisfied
for the Hopf algebras $u(\D,\lambda,\mu)$.

In the next definition we extend the notation for the linking
parameters by \eqref{l2}.

\begin{Def}\label{DefIso}
Let $\G,\G'$ be abelian groups and let
\begin{align*}
\D &= \D(\Gamma, (g_i)_{1 \leq i \leq \theta}, (\chi_i)_{1
\leq i \leq \theta}, (a_{ij})_{1\leq i,j \leq \theta}),\\
\D' &=
\D(\Gamma', (g'_i)_{1 \leq i \leq \theta'}, (\chi'_i)_{1 \leq i
\leq \theta'}, (a'_{ij})_{1\leq i,j \leq \theta'})
\end{align*} be Cartan data
of finite type satisfying \eqref{orderqij} and \eqref{orderNJ}. Assume $\theta = \theta'$.
Let $\lambda$ and $\lambda'$ be linking parameters, and $\mu$ and
$\mu'$ root vector parameters for $\D$ and $\D'.$

Let $\varphi : \G' \to \G$ be a group isomorphism, $\sigma \in
S_{\theta}$  a permutation and $(s_i)_{1 \leq i \leq \theta}$ a
family of non-zero elements in $k$. The triple $(\varphi, \sigma,
(s_i))$ is called an {\em isomorphism} from $(\D',\lambda',\mu')$
to $(\D,\lambda,\mu)$ if the following five conditions are
satisfied:
\begin{align}
\varphi(g'_i)&=g_{\sigma(i)} \text{ for all } 1 \leq i \leq \theta.\label{I1}\\
\chi'_i&=\chi_{\sigma(i)}\varphi \text{ for all } 1 \leq i \leq \theta.\label{I2}\\
a'_{ij} &=a_{\sigma(i),\sigma(j)} \text{ for all } 1 \leq i,j \leq \theta.\label{I3}\\
\lambda'_{ij} &= s_i s_j \lambda_{\sigma(i)\sigma(j)} \text{ for all } 1 \leq i,j \leq \theta, i \nsim j\label{I4}
\end{align}

To formulate the fifth condition we have to introduce more
notations for the connected components $J$ of $\D$.
\begin{enumerate}
\item Let $q_J = (\chi_j(g_i))_{i,j \in J}$ be the braiding matrix of
the restriction $\D_J$ of $\D$ to $J$.
\item Using \eqref{I3} we identify the root systems of $J$ and of
$\sigma^{-1}(J)$ with the corresponding root system of the
standard Cartan matrix. Then the restriction of $\sigma$ to
$\sigma^{-1}(J)$ becomes a diagram automorphism $\sigma_J$ of the
corresponding standard Cartan matrix.
\item For any $\beta \in \Phi_J^+$ let $u'_{\beta}(\mu')$ and
$u_{\beta}(\mu)$ be the elements in the group algebras $k[\G']$
and $k[\G]$ defined in \eqref{ualphageneral}. For any family
$a=(a_{\beta})_{\beta \in \Phi_J^+}$ of natural numbers $a_{\beta}
\geq 0$ we define the product $u(\mu)^a = \prod_{\beta \in
\Phi_J^+} u_{\beta}(\mu)^{a_{\beta}}$.
\item For any $\alpha \in \Phi_J^+$ with  $\alpha = \sum_{i \in J} n_i
\alpha_i, n_i \geq 0 \text{ for all } 1 \leq i \leq \theta$, let
$s_{\alpha} = \prod_{i \in J} s_{i}^{n_i}$.
\item For any $\alpha \in \Phi_J^+$ and any family $a=(a_{\beta})_{\beta
\in \Phi_J^+}$ of natural numbers $a_{\beta} \geq 0$ let
$t_{\alpha, q_J,\sigma_J}^a$ be the element in $k$ defined below
in Theorem \ref{constants} applied to $\D_J$.
\end{enumerate}

Then the last condition is the following identity in the group
algebra $k[\G]$:
\begin{equation}\label{I5}
\varphi(u'_{\alpha}(\mu')) = s_{\alpha}^{N_J}\sum_a
t_{\alpha, q_J,\sigma_J}^a u(\mu)^a \text{ for all } \alpha \in \Phi_J^+, J \in \mathcal{X}.
\end{equation}
Finally let $\Isom((\D',\lambda',\mu'),(\D,\lambda,\mu))$ be the
set of all isomorphisms from $(\D',\lambda',\mu')$ to
$(\D,\lambda,\mu)$.

For Hopf algebras $A',A$ we denote by $\Isom(A',A)$ the set of all
Hopf algebra isomorphisms from $A'$ to $A$.
\end{Def}

We now can state the main result of this section.

\begin{Thm}\label{iso}
Let $\D$ and $\D'$ be Cartan data of finite type with finite
abelian groups $\G$ and $\G'$ and rank $\theta$ and $\theta'$.
Assume that $\D$ and $\D'$ satisfy \eqref{orderqij},
\eqref{orderNJ}. In addition  assume the following condition on
the braiding matrix of $\D$:
\begin{equation}\label{>4}
\ord(q_{ii}) >4 \text{ for all } 1 \leq i \leq \theta.
\end{equation}
Let $\lambda$ and $\lambda'$ be linking parameters, and $\mu$ and
$\mu'$ root vector parameters for $\D$ and $\D'$.

If the Hopf algebras $u(\D',\lambda',\mu')$ and
$u(\D,\lambda,\mu)$ are isomorphic, then $\theta' = \theta$.
Assume $\theta' = \theta$. Then the map
$$ \Isom((\D',\lambda',\mu'),(\D,\lambda,\mu)) \to
\Isom(u(\D',\lambda',\mu'),u(\D,\lambda,\mu))$$ given by
$(\varphi, \sigma, (s_i)) \mapsto F$, where $F(x'_i) = s_i
x_{\sigma(i)}$ and $F(g') = \varphi(g')$ for all $1 \leq i \leq
\theta$ and $g' \in \G'$, is bijective.
\end{Thm}
Before we begin with the proof of Theorem \ref{iso} we need some preparations.

First we see that condition \eqref{I3} is in most cases redundant.

\begin{Lem}\label{redundant}
In the situation of Definition \ref{DefIso} assume $\theta' =
\theta$, \eqref{I1}, \eqref{I2} and \eqref{>4}. Then \eqref{I3}
holds.
\end{Lem}
\pf For all $1 \leq i,j \leq \theta$ let $q'_{ij} =
\chi'_j(g'_i)$, $q_{ij} =\chi_j(g_i)$. Then for all $i,j$
\eqref{I1} and \eqref{I2} imply that $q'_{ij} =
q_{\sigma(i)\sigma(j)}$. Hence $a_{ij} = a'_{\sigma(i)\sigma(j)}$,
since $q_{ii}^{a_{ij}} = q_{ii}^{a'_{\sigma(i)\sigma(j)}},$ and
$a_{ij} - a'_{\sigma(i)\sigma(j)} \in \{0,\pm1,\pm2,\pm3\}.$ \epf

We need an extra information in the situation of Theorem \ref{R(D)}.

\begin{Lem}\label{kernel}
Let $\D$ be a connected Cartan datum of finite type with root
system $\Phi$, finite abelian group $\G$ and $N = \ord(q_{ii})$
for all $i$. Assume \eqref{orderqij} and \eqref{orderNJ}. By
Theorem \ref{CartanNichols} there is a canonical projection $\pi :
R(\D) \to \mathfrak{B}(V)$ whose kernel is the ideal generated by
all $x_{\alpha}^N, \alpha \in \Phi^+$. We denote the coalgebra
structure of $R(\D)$ by $\Delta(x) = x^{(1)} \otimes x^{(2)}$ for
all $x \in R(\D)$. Then
$$K(\D) = R(\D)^{\co \pi} = \{ x \in R(\D) \mid x^{(1)} \otimes \pi(x^{(2)}) = x \otimes 1 \}.$$
\end{Lem}
\pf This follows by bosonization from the corresponding result for
pointed Hopf algebras in \cite{Ma}. \epf

In the following  theorem we define the constants in \eqref{I5}.
\begin{Thm}\label{constants}
Let $\D(\Gamma, (g_i)_{1 \leq i \leq \theta}, (\chi_i)_{1 \leq i
\leq \theta}, (a_{ij})_{1\leq i,j \leq \theta})$ be a connected
Cartan datum of finite type with root system $\Phi$, finite
abelian group $\G$, and $N = \ord(q_{ii}),1 \leq i \leq \theta$.
Assume \eqref{orderqij} and \eqref{orderNJ}. Assume that $(a_{ij})_{1\leq i,j \leq \theta}$ is a standard
Cartan matrix, and let $\sigma$ be a diagram automorphism of
$(a_{ij})_{1\leq i,j \leq \theta}$. Define
$$\D^{\sigma} = \D(\G,(g'_i)_{1 \leq i \leq \theta}, (\chi'_i)_{1
\leq i \leq \theta},(a_{ij})_{1\leq i,j \leq \theta})$$ with $g'_i
= g_{\sigma(i)}, \chi'_i = \chi_{\sigma(i)}$ for all $1 \leq i
\leq \theta$. Let $V \in {\YDG}$ with basis $x_i \in
V_{g_i}^{\chi_i}$, and $V^{\sigma} \in {\YDG}$ with basis
$x_i^{\sigma} \in (V^{\sigma})_{g'_i}^{\chi'_i}$ for all $1 \leq i
\leq \theta$. Then there is an algebra map
$$F^{\sigma} : R(\D^{\sigma}) \to R(\D), x_i^{\sigma} \mapsto
x_{\sigma(i)} \text{ for all } 1 \leq i \leq \theta.$$
For each $\alpha \in \Phi^+,$ and each family $a =
(a_{\beta})_{\beta \in \Phi^+}$ of natural numbers $a_{\beta} \geq
0$ there are uniquely determined elements $t_{\alpha}^a \in k$
depending on the braiding matrix of $\D$ and the diagram
automorphism $\sigma$ such that
$$F^{\sigma}(x_{\alpha}^{\sigma})^N = \sum_a t_{\alpha}^a z^a,$$
where $z^a = x_{\beta_1}^{a_1N} \cdots x_{\beta_p}^{a_pN}$ as in
Definition \ref{Defza} with the fixed ordering
$\beta_1,\dots,\beta_p$ of the positive roots, and where
$x_{\alpha}^{\sigma}$ denotes the root vector of $\alpha$ in
$R(\D_{\sigma})$.
\end{Thm}
\pf The Cartan condition is satisfied for $\D^{\sigma}$ since
$\sigma$ is a diagram automorphism. The linear map 
$$f^{\sigma} :
V^{\sigma} \to V, x_i^{\sigma} \mapsto x_{\sigma(i)}, 1 \leq i
\leq \theta,$$ 
is an isomorphism of Yetter-Drinfeld modules. Then
$f^{\sigma}$ induces isomorphisms $F^{\sigma} : R(\D^{\sigma}) \to
R(\D)$ and $\overline{F^{\sigma}} : \mathfrak{B}(V^{\sigma}) \to
\mathfrak{B}(V)$. Let $\pi : R(V) \to \mathfrak{B}(V)$ and
$\pi^{\sigma} : R(V^{\sigma}) \to \mathfrak{B}(V^{\sigma})$ be the
natural projections. Since $K(\D) = R(\D)^{\co \pi}$ and
$K(\D^{\sigma}) = R(\D^{\sigma})^{\co \pi^{\sigma}}$ by Lemma
\ref{kernel}, it follows that $F^{\sigma}$ maps $K(\D^{\sigma})$
into $K(\D)$. This proves the claim by Theorem \ref{R(D)}. \epf
The meaning of the elements $F^{\sigma}(x_{\alpha}^{\sigma})$ in
the previous theorem can be explained as follows. Let $x_{\alpha}$
be represented as iterated skew-commutator of simple root vectors
$x_{i_1},\dots,x_{i_s}$ in this order. Then
$F^{\sigma}(x_{\alpha}^{\sigma})$ is the same iterated
skew-commutator of the sequence
$x_{\sigma(i_1)},\dots,x_{\sigma(i_s)}$.

As an example, take the Dynkin diagram $A_2$ with the non-simple
root $\alpha = \alpha_1 + \alpha_2$ and the diagram automorphism
$\sigma$ with $\sigma(1) = 2,\sigma(2) =1$. Then $x_{\alpha} =
x_1x_2 - q_{12}x_2x_1$ and $F^{\sigma}(x_{\alpha}^{\sigma}) =
x_2x_1 - q_{21} x_1x_2$.

Finally we note
\begin{Lem}\label{Lemmadifferent}
Let $\D = \D(\Gamma, (g_i)_{1 \leq i \leq \theta}, (\chi_i)_{1
\leq i \leq \theta}, (a_{ij})_{1\leq i,j \leq \theta})$ be a datum
of finite Cartan type and assume \eqref{>4}. Then for all $1 \leq
i,j \leq \theta, i \neq j$, $g_i \neq g_j \text{ or } \chi_i \neq
\chi_j$.
\end{Lem}
\pf Assume there are $i \neq j$ with $g_i=g_j,\chi_i=\chi_j.$ Then
$q_{ii} = q_{jj},$ and $q_{ii}^2 = q_{ii}^{a_{ij}} =
q_{ii}^{a_{ji}}.$ Hence in contradiction to our assumption we have
$q_{ii}^{2-a_{ij}} =1$ for $a_{ij} \in \{0,-1,-2\}$, and
$q_{ii}^{a_{ij}-a_{ji}}=1$ for $a_{ij}=-3$. \epf

We can now prove Theorem \ref{iso}:

\pf Assume that there is a Hopf algebra isomorphism
$$F : A' =u(\D',\lambda',\mu') \to A=u(\D,\lambda,\mu).$$
Then $F$ preserves the coradical filtration and induces an
isomorphism $A'_0=k[\G'] \cong A_0=k[\G],$ given by a group
isomorphism $\varphi : \G' \to \G,$ and by Corollary \ref{gru} an
isomorphism
$$A'_1= k[\G'] \oplus \bigoplus_{\substack{g' \in \G',\\ 1 \leq i
\leq \theta'}} kx'_ig' \cong A_1 =k[\G] \oplus
\bigoplus_{\substack{g \in \G,\\ 1 \leq i \leq \theta}} kx_ig .$$
Hence it follows from Lemma \ref{Lemmadifferent} (see
\cite[6.3]{AS2}) that $\theta=\theta',$ and that there are a
permutation $\sigma \in S_{\theta}$ and elements $0\neq s_i \in k,
1 \leq i \leq \theta$ such that \eqref{I1} and \eqref{I2} hold,
and $F(x'_i) = s_i x_{\sigma(i)}$ for all $1 \leq i \leq \theta$.
Then Lemma \ref{redundant} implies \eqref{I3}, and
$F([x'_i,x'_j]_{c'}) = s_is_j[x_{\sigma(i)},x_{\sigma(j)}]_{c}$
for all $1 \leq i \leq \theta$. Now \eqref{I4} follows from the
linking relations.

To establish \eqref{I5} we fix a connected component $J$ of $\D$,
and we identify $J$ and $J' = \sigma^{-1}(J)$ with the index set
$\{1,\dots,\theta_J\}$ of the corresponding standard Cartan
matrix. Then the restriction of $\sigma$ to $J'$ becomes the
diagram automorphism $\sigma_J$. Let $V_J \in {\YDG}$ with basis
$x_i \in (V_J)^{\chi_i}_{g_i}, i \in J$, and $V'_{J} \in {\YDg}$
with basis $x'_i \in (V'_{J})^{\chi'_i}_{g'_i}, i \in J$. Let
$f^{\sigma_J}_J : V_J^{\sigma} \to V_J$ be the map of Theorem
\ref{constants} for $\D_J$ and $V_J$ instead of $\D$ and $V$.

We define linear maps $f_J: V'_J \to V_J, x'_i \mapsto
s_ix_{\sigma(i)}, 1 \leq i \leq \theta_J,$ and $f'_J: V'_J \to
V_J^{\sigma}, x'_i \mapsto s_ix^{\sigma_J}_i, 1 \leq i \leq
\theta_J.$ Then $f_J = f^{\sigma_J}_Jf'_J$. The maps $f_J$ and
$f'_J$ are $\G'$-linear and $\G$-colinear, where the action of
$\G'$ on $V_J$ and the coaction of $\G$ on $V'_J$ are defined via
$\varphi$. Hence by \eqref{I1}--\eqref{I3} they induce algebra
maps $F_J :R(\D'_{J'}) \to R(\D_J)$ and $F'_J :R(\D'_{J'}) \to
R(\D_J^{\sigma_J})$, and $F_J = F^{\sigma_J}_JF'_J$. By Theorem
\ref{constants} we then obtain for any $\alpha \in \Phi_J^+$
\begin{equation}\label{final}
F_J(x'_{\alpha})^{N_J}= s_{\alpha}^{N_J}F_J^{\sigma_J}(x^{\sigma_J}_{\alpha})^{N_J} =s_{\alpha}^{N_J}\sum_{a}
t_{\alpha, q_J,\sigma_J}^a z^a.
\end{equation}

The canonical maps
$$\pi_J : R(\D_J) \to u(\D,\lambda,\mu) \text{ and }
\pi_{J'} : R(\D'_{J'}) \to u(\D',\lambda',\mu')$$
map root vectors
to root vectors. Since $F\pi_{J'} = \pi_JF_J$, we see from
\eqref{final} that for all $\alpha \in \Phi_J^+$
$$F\pi_{J'}({x'}_{\alpha}^{N_J}) = \pi_JF_J({x'}_{\alpha}^{N_J})= s_{\alpha}^{N_J}\sum_{a}
t_{\alpha, q_J,\sigma_J}^a u(\mu)^a.$$ On the other hand
${x'}_{\alpha}^{N_J} = {u'}_{\alpha}(\mu')$ in
$u(\D',\lambda',\mu')$, hence
$$F\pi_{J'}({x'}_{\alpha}^{N_J})=\varphi(u'_{\alpha}(\mu')),$$ and
\eqref{I5} follows.

It is easy to see that conversely any isomorphism
$(\varphi,\sigma,(s_i))$ defines an isomorphism of Hopf algebras,
and that two such triples coincide if they define the same Hopf
algebra map. \epf We remark that the situation greatly simplifies
if the diagram automorphism $\sigma_J$  in Definition \ref{DefIso}
is the identity. This happens in particular if the Dynkin diagram
of $(a_{ij})_{i,j \in J}$ is not of Type $A, D$ or $E_6$. In this
case it follows from the inductive definition of the
$u_{\alpha}(\mu)$ that \eqref{I5} is equivalent to
\begin{equation}
\mu'_{\alpha} =s_{\alpha}^{N_J} \mu_{\alpha} \text{ for all }\alpha \in \Phi_J^+,J \in \mathfrak{X}.
\end{equation}

\end{document}